%% file: main_elsevier.tex
\documentclass[authoryear,preprint,11pt,3p]{elsarticle}
\journal{Computers \& Chemical Engineering}

\usepackage{graphicx}
\usepackage{xcolor}
\usepackage{physics}
\usepackage{amsmath}
\usepackage{amssymb}
\usepackage{amsthm}
\usepackage{eurosym}
\theoremstyle{plain}%
\newtheorem{theorem}{Theorem}
\newtheorem{lemma}{Lemma}
\theoremstyle{definition}

\usepackage{ulem}

\usepackage{nomencl}
\makenomenclature
\usepackage{longtable}
\makeatletter
\def\els@aparagraph[#1]#2{\elsparagraph[#1]{#2\@addpunct{.}}}
\def\els@bparagraph#1{\elsparagraph*{#1\@addpunct{.}}}
\makeatother
\usepackage{bm}
\usepackage{amsfonts}
\usepackage[colorlinks=True]{hyperref}
\usepackage{booktabs}
\usepackage{multirow}
\usepackage{mathtools}
\usepackage[bb=boondox]{mathalfa}
\usepackage{cleveref}
\crefformat{section}{§#2#1#3}
\usepackage{nicefrac}
\usepackage{enumitem}
\usepackage{chngcntr}
\usepackage{booktabs}
\usepackage{tikz,tikz-3dplot}
\usepackage{pgfplots,pgfplotstable}
\usetikzlibrary{arrows,automata,calc,shapes,plotmarks, positioning,shapes.geometric,decorations.pathmorphing,patterns,fillbetween,pgfplots.groupplots}
\usepackage{circuitikz}
\pgfplotsset{compat=1.15}
\pgfplotsset{
    table/search path={plotdata},
}
\usepgfplotslibrary{colorbrewer,units}
\pgfplotsset{cycle list/Set3}

\definecolor{maincolor}{HTML}{032F99} 
\definecolor{secondcolor}{HTML}{ff5722} 
\definecolor{thirdcolor}{HTML}{c7d3d7}  
\definecolor{newtextcolor}{HTML}{bf360c}

\usepackage{adjustbox}

\usepackage{mathtools}

\usepackage{esvect}

\makeatletter
\newcommand{\oset}[3][0ex]{%
  \mathrel{\mathop{#3}\limits^{
    \vbox to#1{\kern-2\ex@
    \hbox{$\scriptstyle#2$}\vss}}}}
\makeatother

\begin{document}

\begin{frontmatter}
\title{A Conic Model for Electrolyzer Scheduling}
\author[inst1]{Enrica Raheli\corref{cor1}}\ead{enrah@dtu.dk}
\author[inst1,inst2]{Yannick Werner\corref{cor1}}\ead{yanwe@dtu.dk}
\author[inst1]{Jalal Kazempour}\ead{jalal@dtu.dk}

\affiliation[inst1]{organization={Technical University of Denmark},
                    city={Kgs. Lyngby},
                    Denmark}
\affiliation[inst2]{organization={Norwegian University of Science and Technology}, city={Trondheim}, country={Norway}}
\cortext[cor1]{Corresponding authors. The first two co-authors contributed equally and are listed alphabetically.}

\begin{abstract}
The hydrogen production curve of the electrolyzer describes the nonlinear and nonconvex relationship between its power consumption and hydrogen production. An accurate representation of this curve is essential for the optimal scheduling of the electrolyzer. The current state-of-the-art approach is based on piecewise linear approximation, which requires binary variables and does not scale well for large-scale problems. To overcome this barrier, we propose two models, both built upon convex relaxations of the hydrogen production curve. The first one is a linear relaxation of the piecewise linear approximation, while the second one is a conic relaxation of a quadratic approximation. Both relaxations are exact under prevalent operating conditions. We prove this mathematically for the conic relaxation. Using a realistic case study, we show that the conic model, in comparison to the other models, provides a satisfactory trade-off between computational complexity and solution accuracy for large-scale problems. 
\end{abstract}




\end{frontmatter}
\input{00_Nomenclature.tex}
\input{01_introduction.tex}
\input{02_Preliminaries.tex}

\input{03_Model.tex}
\input{03A_Tightness.tex}
\input{04_Results.tex}
\input{05_Conclusion.tex}
\appendix
\input{Appendix_proofs.tex}
\input{Acknowledgment.tex}

\bibliographystyle{elsarticle-harv} 
\bibliography{references}

\end{document}

%% file: 00_Nomenclature.tex
\section*{Nomenclature}
\vspace{-4mm}
{\small
 \begin{longtable}{p{3.1cm} p{12.5cm}}\label{tab:nomenclature}\\
		\multicolumn{2}{l}{\textbf{Indices and sets}}  \\
		$t \in \mathcal{T}$         & Set of time steps $t$ \\
        $n \in \mathcal{N}$         & Set of subperiods $n$ \\
        $\mathcal{H}_n \subseteq \mathcal{T}$         & Set of time steps in subperiod $\textit{n}$ \\
        $s \in\mathcal{S}$         & Set of linearization segments $s$ \\
        & \\
        \multicolumn{2}{l}{\textbf{Parameters}}  \\
        $\lambda_t$ & Day-ahead power price in time step~$t$ [\euro/MWh]\\
        $\chi$ & Hydrogen price [\euro/kg]\\
        $A_s, B_s, C_s$  & Polynomial coefficients for segment $s$ $\rm{[kg/(hMW^2)],[kg/(hMW)],[kg/h]}$ \\
        $D^{\mathrm{max}}_{n}$ & Maximum hydrogen production in subperiod $n$ [kg]\\
        $K^{\mathrm{su}}$ & Cold startup cost of the electrolyzer [\euro]\\
        $P^{\mathrm{min}}, P^{\mathrm{max}}, P^{\mathrm{sb}}$ & Minimum, maximum, and standby power consumption of the electrolyzer [MW]\\
        $\underline{P}_s, \overline{P}_s $ & Lower and upper bounds for power consumption of segment $s$ [MW]\\
        $W_t$ & Day-ahead wind power production forecast in time step~$t$ [MW]\\
        $Q_1, Q_0$  & Polynomial coefficients of the linear underestimator $\mathrm{[kg/(hMW)],[kg/h]}$ \\
        & \\
        \multicolumn{2}{l}{\textbf{Variables}}  \\
        $f_t \in \mathbb{R^+}$ & Power sold to the grid in time step~$t$ [MW]\\
        $h_t \in \mathbb{R^+}$ & Hydrogen production of the electrolyzer in time step~$t$ [kg/h]\\
        $p_t \in \mathbb{R^+}$ & Power consumption of the electrolyzer in time step~$t$ [MW]\\
        $\widehat{z}_{s,t} \in \{0,1\}$ & Binary variable defining the segment $s$ that the electrolyzer operates on in time step~$t$\\ 
        $\widehat{p}_{s,t}\in \mathbb{R^+}$ & Power consumption of the electrolyzer in on state in segment $s$ and time step~$t$ [MW]\\
        $\widetilde{p}_{t} \in \mathbb{R^+}$ & Power consumption of the electrolyzer in on state in time step~$t$ [MW]\\
        $z^{\mathrm{on}}_{t}, z^{\mathrm{off}}_{t}, z^{\mathrm{sb}}_{t} \in \{0,1\}$ & Binary variables defining on, off, and standby state of the electrolyzer in time step~$t$ \\
        $z^{\mathrm{su}}_{t} \in \{0,1\}$ & Binary variable defining a cold startup of the electrolyzer in time step~$t$ \\ 
	\end{longtable}
}

\subsection*{Notes on the nomenclature}
Variables (continuous and binary) are denoted by lower-case letters. In particular, binary variables are denoted by the letter $z$.
Parameters are denoted by upper-case letters, except hydrogen and electricity prices which are denoted by Greek letters.
The set $\mathbb{R^+}$ denotes the set of nonnegative real numbers.

%% file: 01_introduction.tex
\section{Introduction}
\label{sec:intro}
\subsection{Background}
Renewable hydrogen produced via electrolysis is widely acknowledged as a key priority to achieve a clean energy transition. Several countries in Europe and globally have published national hydrogen strategies to support the large-scale development of electrolyzers \citep{NationalHydrogenStategies}. For instance, the 2020 EU Hydrogen strategy sets an electrolyzer capacity target of 40 GW by 2030 \citep{PtX_strategy_EU}. 
Developing this technology on a large scale poses several challenges, including the scale-up of manufacturing processes, the improvement of system design and materials, the establishment of a supportive policy framework, and the definition of new business models \citep{cost_reduction, newbusinessmodels}.

Hybrid power plants consisting of renewable power sources (wind and/or solar) and electrolyzers create synergies based on cross-commodity arbitrage between electricity and hydrogen markets \citep{newbusinessmodels}.
The plant operator can dynamically control the system to take advantage of volatile electricity prices: selling electricity directly when power prices are high and producing and selling hydrogen when power prices are low \citep{MATUTE20211449}. In this way, the cost of hydrogen production, which mainly depends on the cost of electricity \citep{HenrikLundFrandsen}, is reduced.
This requires the development of optimal scheduling models aiming at maximizing the profit of the hybrid power plant.

To accurately capture the operational space of the hybrid power plant, those models should be aware of the underlying physics of its components. This poses a challenge particularly for modeling the electrolyzer, for which the relationship between its power consumption and hydrogen production is nonlinear and nonconvex. We refer to this relationship as the \textit{hydrogen production curve}, which does not have a known analytical expression.\footnote{In \cite{Ulleberg}, a widely adopted empirical equation is presented to describe the relationship between electrolyzer voltage and current density. This equation is noninvertible, making it impossible to derive an analytical expression for the hydrogen production as a function of the power consumption of the electrolyzer. For a comprehensive explanation of the methodology used in this paper to construct the hydrogen production curve, we direct the reader to \cite{powertech2023}.} 
As the technology is still in an early stage of development, manufacturers disclose very limited information on technical characteristics. From modeling perspective, two main questions arise: First, how to best approximate the hydrogen production curve under limited information availability? Second, how to deal with computational complexity of the nonlinear hydrogen production curve when optimizing dispatch strategies?
\subsection{Electrolyzer Hydrogen Production Modeling: Status Quo} 
\begin{figure}[t]
\centering
\includegraphics[width=0.5\linewidth]{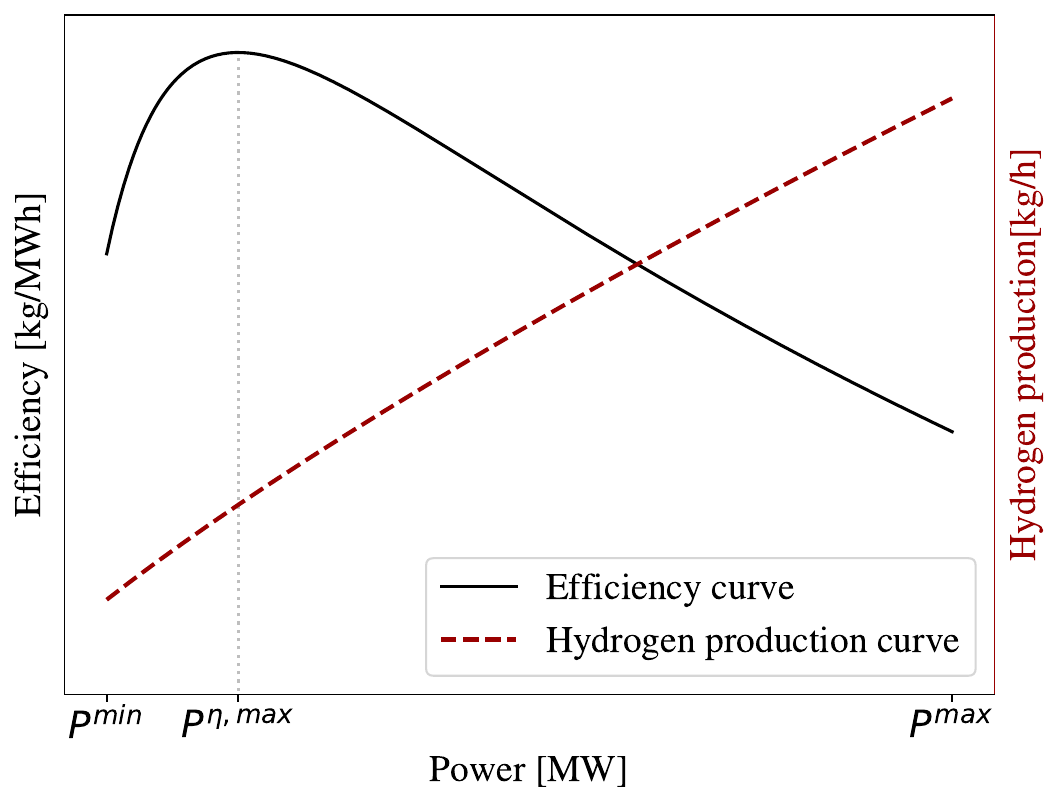}%
\caption{Schematic efficiency (black) and hydrogen production (red) curves for an alkaline electrolyzer. $P^{\mathrm{min}}$ and $P^{\mathrm{max}}$ refer to the minimum and maximum operating range, respectively. $P^{\mathrm{\eta, max}}$ is the power consumption corresponding to the peak efficiency. This figure is plotted based on the methodology in \cite{powertech2023} applied to data from \cite{Ulleberg} and \cite{Sanchez}.} 
\label{fig:curves}
\end{figure}

The dashed red curve in Figure~\ref{fig:curves} shows a schematic hydrogen production curve of an alkaline electrolyzer. To better illustrate the nonlinear physics of the electrolyzer, the black curve depicts the ratio of hydrogen production to power consumption, the so-called \textit{efficiency curve}.
The efficiency peaks at around $20$-$40$\% of the electrolyzer capacity \citep{Siemens_white_paper}, for which the corresponding power consumption of the electrolyzer is denoted as $P^{\mathrm{\eta, max}}$. For power consumption levels higher than $P^{\mathrm{\eta, max}}$, the efficiency declines almost linearly due to the effect of overpotentials in the polarization curve \citep{SANCHEZ20203916}. For consumption levels below $P^{\mathrm{\eta, max}}$, the efficiency drops rapidly due to increasing Faraday losses \citep{Ulleberg}. This nonlinearity needs to be incorporated into operational decision-making problems of the electrolyzer through the hydrogen production curve (red). For a more detailed analysis of the physics and technical modeling of the electrolyzer, the interested reader is referred to \citet{Ulleberg}, \citet{Sanchez}, and \citet{ShiYou}.

Currently, there is no widely adopted approach for modeling the hydrogen production curve in a computationally tractable way. It is a common practice in the literature to introduce relatively strong simplifications for the hydrogen production curve. For instance, a constant efficiency is used in \cite{MATUTE20211449} and \cite{pavic2022}. In \cite{VARELA20219303}, a first and second-order polynomial approximation for the hydrogen production curve is proposed. Although the second-order polynomial exhibits a smaller error, it has been eventually discarded due to its computational complexity. In \citet{powertech2023}, it has been shown that a linear hydrogen production curve is not suitable when the electrolyzer frequently operates at partial load, i.e., in the range between $P^{\mathrm{min}}$ and $P^{\mathrm{max}}$ in Figure~\ref{fig:curves}. Instead, they propose a piecewise linear hydrogen production curve based on an incremental method and assess the impact of neglecting a detailed electrolyzer model. However, additional binary variables, one for each linearization segment, need to be introduced, increasing the computational burden. Similar piecewise linearization approaches based on convex combination methods are used in \citet{YiZeng_PWL}, \citet{Lin2021}, and \citet{Kountouris2023}. The reader is referred to \citet{PWLmethods} to learn more about different methods for the piecewise linear approximation.

Convex relaxation techniques of nonlinear physics are often employed in the energy system literature. For instance, extensive research has been conducted on developing convex relaxations of the AC optimal power flow problem \citep{ZOHRIZADEH202020} and for the optimal natural gas flow problem \citep{Pascal2016}. To the best of our knowledge, convex relaxation techniques for modeling the nonlinear hydrogen production curve of electrolyzers have not been explored.

\subsection{Contributions and Paper Organization}
Existing approaches to modeling the nonlinear hydrogen production curve are either inaccurate, resulting in sub-optimal operational decisions, or do not computationally scale well for large optimization problems. Some examples of large-scale problems involving electrolyzers are the operation and planning of complex power-to-x systems with further downstream chemical processes; the trading problem in multiple markets for electricity, hydrogen, and ancillary services; the operation of individual stacks of large-scale electrolyzers; and the scheduling problem under uncertainty, e.g., when the electrolyzer is part of a renewable-based hybrid power plant.

This paper takes a new perspective on modeling the hydrogen production curve based on convex relaxations, ensuring both accuracy and computational scalability for large-scale problems.
Starting from the state-of-the-art piecewise linearization method, we derive a corresponding linear relaxation that does not require binary variables.
As our main contribution, we propose a conic relaxation of a quadratic approximation of the hydrogen production curve, resulting in the so-called \textit{conic model}.
We mathematically define the sufficient conditions for the conic relaxation to be exact and the necessary and sufficient conditions for it to be inexact. 
We demonstrate those findings using an illustrative case study.
The conic model can be readily applied to the large-scale optimization problems mentioned above, among others. We numerically compare it to the state-of-the-art piecewise linearization method and its linear relaxation counterpart, and verify ex-post its well-performance in terms of solution quality, i.e., operational decisions, and reduced computational complexity.

The remainder of this paper is organized as follows. Section~\ref{sec:prelim} introduces the system setup and the optimal scheduling problem of the hybrid power plant. Section~\ref{sec:model} provides an overview of the different approaches to model the nonlinear hydrogen production curve, including the proposed relaxations. In Section~\ref{sec:tightness}, we prove the exactness of the proposed conic relaxation. Section~\ref{sec:results} compares the modeling approaches, and Section~\ref{sec:conclusion} concludes.

%% file: 02_Preliminaries.tex
\section{Scheduling Problem}
\label{sec:prelim}

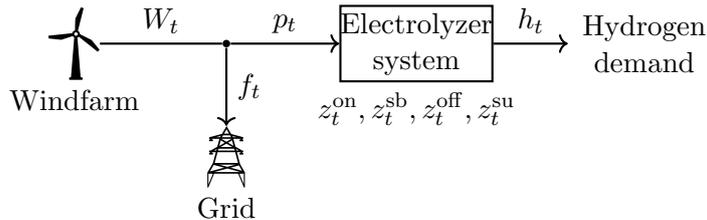
\begin{figure}[t]
\centering
\begin{tikzpicture}
\node[inner sep=0pt] (wind) at (0,0) [label=below:{Windfarm}]
    {\includegraphics[width=.04\textwidth]{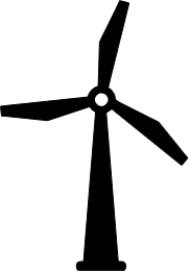}};
\node[inner sep=0pt] (grid) at (2,-1.5) [label=below:{Grid}]
    {\includegraphics[width=.03\textwidth]{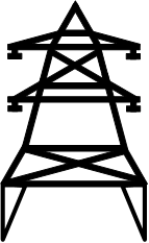}};
\node[inner sep=0pt] (electrolyzer) at (4.5,0)
   [draw,thick,align=center,minimum width=2cm,minimum height=1cm, text width = 2 cm, label=below:{$z^{\mathrm{on}}_{t},z^{\mathrm{sb}}_{t},z^{\mathrm{off}}_{t},z^{\mathrm{su}}_{t}$}] {Electrolyzer system};
\node[inner sep=0pt,fill,circle,minimum size=0.1 cm] (meter) at (2,0) {};
\node[inner sep=0pt] (demand) at (7.5,0) [text width = 2 cm,align=center] {Hydrogen demand};

\draw[-,thick] (wind.east) -- (meter.west)
    node[midway,above] {$W_t$};
\draw[->,thick] (meter.east) -- (electrolyzer.west)
    node[midway,above] {$p_t$};
\draw[->,thick] (electrolyzer.east) -- (demand.west)
    node[midway,above] {$h_t$};
\draw[->,thick] (meter.south) -- (grid.north)
    node[midway,right] {$f_t$};    
\end{tikzpicture}
\caption{Schematic representation of a hybrid power plant including the decision variables of the corresponding optimization problem. The electrolyzer system includes the hydrogen compressor and other auxiliary components.}
\label{fig:system}
\end{figure}

We consider a hybrid power plant as shown in Figure~\ref{fig:system}, consisting of a wind farm, an electrolyzer system, and a hydrogen demand. Hereafter, we refer to the electrolyzer as a system including all necessary auxiliary components, such as transformers, rectifiers, pumps, coolers, heaters, and compressors, in addition to the electrolyzer device. Following the current common practice, we assume that hydrogen is sold at a fixed price and that there is an upper limit on the hydrogen demand for a given time interval. The hydrogen demand could be constrained by storage availability or maximum capacity of downstream facilities, such as industry or chemical plants for further conversion into e-fuels. Here, for simplicity, we assume that the demand is constrained by a tube trailer with a constant capacity and scheduled (e.g., daily) pickups. As we do not consider any constraints on the total hydrogen production directly, it is effectively restricted by the total hydrogen demand. We, therefore, use the terms total hydrogen production and demand interchangeably. For the sake of simplicity, we assume that there is no minimum level of hydrogen demand. Electricity generated by the wind farm can either be sold to the grid at a perfectly known price or used for hydrogen production through the electrolyzer. Since the primary focus of this paper is on the modeling of the hydrogen production curve for electrolyzer scheduling in general, ancillary services or other electricity markets are not considered. Participation in these markets can increase the profitability of the hybrid power plant \citep{Zheng2023,saretta2023}. Without loss of generality and to avoid discussion about the carbon intensity of the hydrogen produced, we assume that the hybrid power plant never buys electricity from the grid.\footnote{According to the ``Renewable Fuels of Non-Biological Origin (RFNBO)" Delegated Act, published by the European Commission in 2023, there are some exceptions that allow the hydrogen production to be labeled green when buying electricity from the grid, e.g., when the price or carbon intensity of the electricity is low. These exceptions are not considered in this paper but the interested reader is referred to \cite{RFNBO} to learn more about the EU regulation.}

In the following, we introduce a model for the optimal day-ahead scheduling problem of the given hybrid power plant. We illustrate the problem by showing the corresponding decision variables in Figure~\ref{fig:system}. The modeling of the hydrogen production curve and the proposed relaxations will be presented in Section~\ref{sec:model}.
To ease notational clarity, we use upper-case and Greek symbols for parameters, and lower-case symbols for variables. Let $t \in \mathcal{T}$ denote the set of time steps, which is divided into $n \in \mathcal{N}$ subsets $\mathcal{H}_{n} \subseteq \mathcal{T}$, such that $\cup_{n \in \mathcal{N}} \mathcal{H}_{n} = \mathcal{T}$ and $\cap_{n \in \mathcal{N}} \mathcal{H}_{n} = \emptyset$. For example, $t \in \mathcal{T}$ could be $8760$ hours of the year, whereas $\mathcal{H}_{n}, \ \forall{n} \in \{1, 2, ..., 365\}$, indicates the set of $365$ days. 

The operator of the hybrid power plant maximizes the total profit from selling power $f_t$ to the grid at the day-ahead power price $\lambda_t$, and selling hydrogen $h_{t}$ at a constant price $\chi > 0$. The only operational expense considered is the startup cost $K^{\mathrm{su}}$ indicated by the binary variable $z^{\mathrm{su}}_{t}$:
\begin{align}
\begin{split} \label{eq:DA_objective}
\underset{\mathbf{x},\mathbf{y},\mathbf{z}}{\mathrm{max}} & \sum_{t \in \mathcal{T}} \left( f_t \lambda_t + h_t \chi  - z^{\mathrm{su}}_{t} K^{\mathrm{su}} \right),
\end{split}
\end{align}
where $\mathbf{x}$, $\mathbf{y}$, and $\mathbf{z}$ denote the set of variables corresponding to the balance of the hybrid power plant, the hydrogen production curve, and the operational states of the electrolyzer, respectively. These three sets will be defined later. 
The power balance within the hybrid power plant is enforced by
\begin{subequations} 
\begin{align}
    & W_t - f_t - p_t = 0, & \forall\,t \in \mathcal{T}, \label{eq:DA_balance}\\
    & f_t \geq 0, & \forall\,t \in \mathcal{T} \label{eq:DA_green_hydrogen},
\end{align}
where $W_t$ denotes the deterministic day-ahead forecast of the wind power production and $p_{t}$ is the day-ahead schedule for the power consumption of the electrolyzer in time step $t$.
Constraint~\eqref{eq:DA_green_hydrogen} prohibits purchasing power from the grid. Without loss of generality, the curtailment of wind power production is not allowed, irrespective of the power price.

In each subset of time steps $\mathcal{H}_{n}$, e.g., over every day $n \in \mathcal{N}$, the total hydrogen production is limited by the hydrogen demand $D^{\mathrm{max}}_{n}$, which, e.g., represents the capacity of an available tube trailer for hydrogen transportation:
\begin{align}
    & \sum_{t \in \mathcal{H}_{n}} h_{t} \leq D^{\mathrm{max}}_{n}, & \forall\, & n \in \mathcal{N}, \label{eq:DA_hydrogen_demand} \\
    & h_{t} \geq 0, & \forall\, & t \in \mathcal{T}. \label{eq:DA_non_negative_hydrogen}
\end{align}
\end{subequations} 
Constraint~\eqref{eq:DA_non_negative_hydrogen} ensures that the hydrogen production is nonnegative.
The set of variables $\mathbf{x}$ is defined as $\mathbf{x} = \{f_t, p_t, h_t\}$.

We consider three operational states for the electrolyzer, namely on, standby, and off \citep{MATUTE20211449, VARELA20219303, ShiYou, powertech2023}. In the on state, the electrolyzer is consuming power and producing hydrogen. Below a certain minimum power consumption, the electrolyzer has to be turned to standby or off. In standby, the electrolyzer does not produce hydrogen but consumes a small amount of power to keep the system running and be able to turn on immediately. On the contrary, in the off state, the electrolyzer does not consume any power but takes several minutes and a significant amount of electricity to be switched back to the on state. The transition from the off to on state is referred to as cold startup. Modeling these operational states and transitions in the scheduling problem requires binary variables. The three operational states, i.e., on, off, and standby, are indicated by binary variables $z^{\mathrm{on}}_{t}, z^{\mathrm{off}}_{t}, z^{\mathrm{sb}}_{t}$, respectively. They constrain the feasible power consumption of the electrolyzer as
\begin{subequations}
\begin{align}
    & z^{\mathrm{on}}_{t} + z^{\mathrm{off}}_{t} + z^{\mathrm{sb}}_{t} = 1, & \forall\, & t \in \mathcal{T}, \label{eq:DA_mutual_states} \\
    & p_{t} \leq P^{\mathrm{max}} z^{\mathrm{on}}_{t} + P^{\mathrm{sb}} z^{\mathrm{sb}}_{t}, & \forall\, & t \in \mathcal{T}, \label{eq:DA_consumption_max} \\
    & p_t \geq P^{\mathrm{min}} z^{\mathrm{on}}_{t} + P^{\mathrm{sb}} z^{\mathrm{sb}}_{t}, & \forall\, & t \in \mathcal{T}, \label{eq:DA_consumption_min} \\
    & z^{\mathrm{su}}_{t} \geq z^{\mathrm{off}}_{t-1} + z^{\mathrm{on}}_{t}  + z^{\mathrm{sb}}_{t} - 1, & \forall\, & t \in \mathcal{T} \backslash 1, \label{eq:DA_start_up} \\
    & z^{\mathrm{su}}_{t=1} = 0, & & \label{eq:DA_start_up_1} \\
    & z^{\mathrm{on}}_{t}, z^{\mathrm{off}}_{t}, z^{\mathrm{sb}}_{t}, z^{\mathrm{su}}_{t} \in \{0,1\} , & \forall\, & t \in  \mathcal{T}. \label{eq:DA_state_binaries}
\end{align}
\end{subequations}
Constraint~\eqref{eq:DA_mutual_states} ensures mutual exclusiveness of the operational states. The corresponding power consumption is constrained by~\eqref{eq:DA_consumption_max} and \eqref{eq:DA_consumption_min} based on the standby, minimum, and maximum power consumption levels $P^{\mathrm{sb}}, P^{\mathrm{min}},$ and $P^{\mathrm{max}}$, respectively. Constraints~\eqref{eq:DA_start_up} and \eqref{eq:DA_start_up_1} define a cold startup $z^{\mathrm{su}}_{t} \in \{0,1\}$ when the electrolyzer changes from off to on state. Lastly, the operational states and transitions are restricted to be binary by~\eqref{eq:DA_state_binaries}. The set of variables $\mathbf{z}$ is defined as $\mathbf{z} = \{z^{\mathrm{on}}_{t}, z^{\mathrm{sb}}_{t}, z^{\mathrm{off}}_{t}, z^{\mathrm{su}}_{t}\}$.

The hydrogen production curve ($\rm{HYP}$), illustrated by the red curve in Figure~\ref{fig:curves}, relates the power consumption to the hydrogen production of the electrolyzer in a general form of
\begin{align}
    & g(h_t,p_t) = 0, & \forall\, t \in \mathcal{T}. \label{eq:DA_hydrogen_general}
\end{align}
In the following section, we present three different approximation and/or relaxation models for \eqref{eq:DA_hydrogen_general}, including our proposed conic model. Each model ends up in a set of constraints that replaces \eqref{eq:DA_hydrogen_general}, which is part of the optimal scheduling problem of the hybrid power plant \eqref{eq:DA_objective}-\eqref{eq:DA_hydrogen_general}.

%% file: 03_Model.tex
\section{Modeling the Hydrogen Production Curve}
\label{sec:model}
There is currently limited information available on the technical characteristics of electrolyzers. Therefore, we compute a hydrogen production curve using the process explained in \cite{powertech2023}, which is based on empirical relationships found by \cite{Ulleberg} and \cite{Sanchez}. We refer to it as the \textit{experimental} hydrogen production curve $\mathrm{HYP}$-$\mathrm{X}$, which does not have a closed-form analytical expression. Note that it is unnecessary to compute $\mathrm{HYP}$-$\mathrm{X}$ if operational data on power consumption and hydrogen production is available, e.g., as in \cite{Energiepark_Mainz}. The models introduced next can then directly be constructed from the operational data.
The three\footnote{In Section~OC.1 of the online companion \citep{Online_Companion}, we further present a fourth model, $\mathrm{HYP}$-$\mathrm{MISOC}$, which has multiple second-order cone constraints.} models that replace \eqref{eq:DA_hydrogen_general} are
\begin{enumerate}
    \item $\mathrm{HYP}$-$\mathrm{MIL}$: The current state-of-the-art piecewise linear approximation,
    \item $\mathrm{HYP}$-$\mathrm{L}$: A corresponding linear relaxation,
    \item $\mathrm{HYP}$-$\mathrm{SOC}$: A second-order cone relaxation of a quadratic approximation.
\end{enumerate}

\begin{figure*}[h]
\centering
\includegraphics[width=\linewidth]{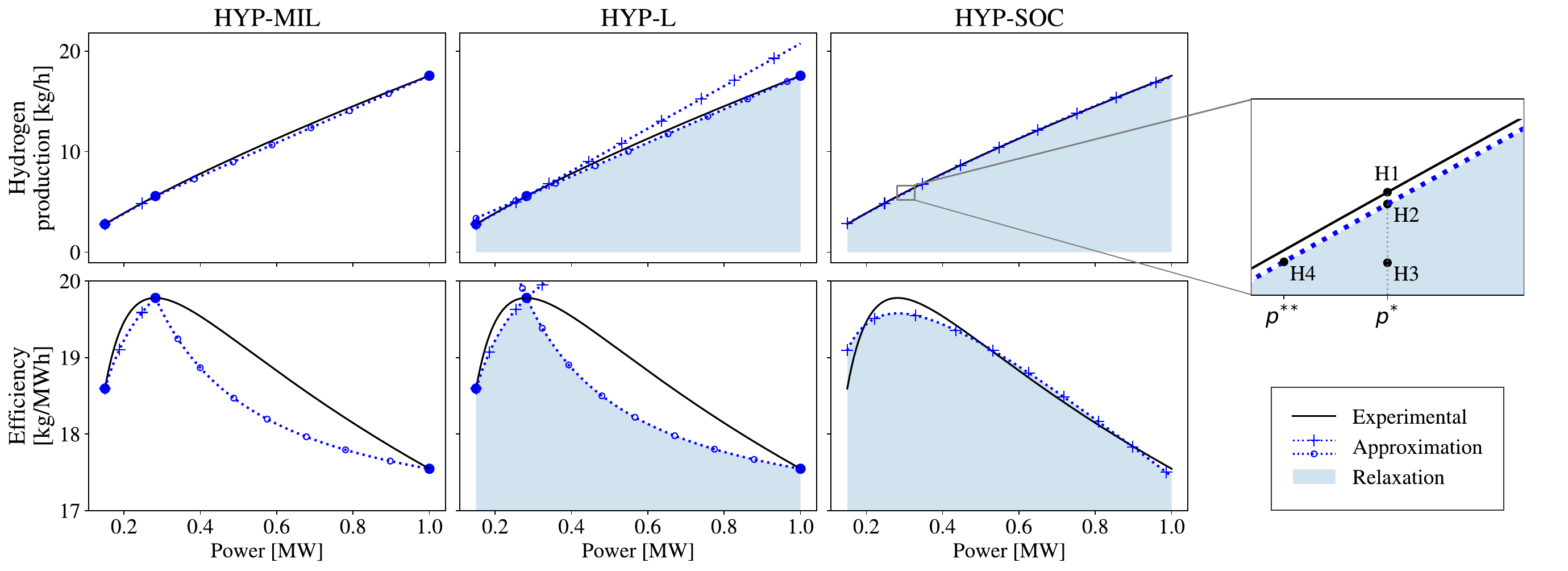}%
\caption{Approximation (blue dashed) and relaxation (blue shaded) of the experimental nonconvex (black) hydrogen production and efficiency curves using two segments (except $\mathrm{HYP}$-$\mathrm{SOC}$). The experimental hydrogen production and efficiency curves are taken from \citet{powertech2023}, based on \citet{Ulleberg} and \citet{Sanchez}. For an exemplary optimal power consumption $p^*$, the actual amount of hydrogen produced according to the experimental curve is indicated by point $\mathrm{H1}$. Points $\mathrm{H2}$ and $\mathrm{H3}$ denote solutions when the relaxation $\rm{HYP}$-$\rm{SOC}$ is exact and inexact, respectively. Point $\mathrm{H4}$ marks the only exact solution that has the same hydrogen production as point $\mathrm{H3}$ but lower power consumption $p^{**}$.} 
\label{fig:3models}
\end{figure*}

These three models along with $\mathrm{HYP}$-$\mathrm{X}$ are illustrated in Figure~\ref{fig:3models}, together with their efficiency curves. The latter is only illustrated to highlight the differences between the approximations and relaxations of the hydrogen production curve but it is not used in the optimization problem.
For illustration clarity, Figure~\ref{fig:3models} shows two segments only for $\mathrm{HYP}$-$\mathrm{MIL}$ and $\mathrm{HYP}$-$\mathrm{L}$, which are indicated by two different dashed lines. In general, any finite number and position of segments can be chosen depending on the desired trade-off between accuracy and computational complexity.
Model $\mathrm{HYP}$-$\mathrm{SOC}$ is a relaxation of a quadratic approximation of the experimental hydrogen production curve such that the corresponding efficiency reaches its maximum at the power consumption $P^{\eta,\mathrm{max}}$.\footnote{We estimate the coefficients of the quadratic approximation by minimizing the sum of squared residuals while giving a relatively higher weight to the residual corresponding to $P^{\eta,\mathrm{max}}$. We found this approximation to deliver better decisions in the scheduling model than just minimizing the sum of squared residuals.} A close-up on model $\mathrm{HYP}$-$\mathrm{SOC}$ is shown in the inset plot on the right side of Figure~\ref{fig:3models}. Let $p^*$ and $h^*$ denote the optimal power consumption and optimal hydrogen production of the electrolyzer, respectively. The actual hydrogen production corresponding to $p^*$, according to the experimental hydrogen production curve $\mathrm{HYP}$-$\mathrm{X}$, is indicated by point $\mathrm{H1}$. Model $\mathrm{HYP}$-$\mathrm{SOC}$ may not be able to attain this point as it utilizes an approximation of the experimental curve, which is illustrated by the blue dashed line. The solution to $\mathrm{HYP}$-$\mathrm{SOC}$ may therefore be located at point $\mathrm{H2}$. We call the discrepancy between $\mathrm{H1}$ and $\mathrm{H2}$ \textit{approximation error}. For models $\mathrm{HYP}$-$\mathrm{MIL}$ and $\mathrm{HYP}$-$\mathrm{L}$, this error is always nonnegative and its magnitude can be controlled by the number of linearization segments. This is not the case for $\mathrm{HYP}$-$\mathrm{SOC}$, where the approximation error can be either positive or negative, and is limited by the shape of the quadratic approximation.
Depending on the choice of the model for the hydrogen production curve, the magnitude of the approximation error is below 1 kg/h.

In contrast to the nonconvex approximation $\mathrm{HYP}$-$\mathrm{MIL}$, models $\mathrm{HYP}$-$\mathrm{L}$ and $\mathrm{HYP}$-$\mathrm{SOC}$ further admit the point $\mathrm{H3}$, located in the blue shaded area, when the relaxation of the hydrogen production curve is inexact. We refer to the discrepancy between $\mathrm{H2}$ and $\mathrm{H3}$ as \textit{relaxation gap}.
The theoretical magnitude of the relaxation gap is upper-bounded  by the hydrogen production curve. When the electrolyzer operates at full capacity, this is around 17.55 kg/h.
The actual magnitude of the gap depends on the exactness of the relaxation, which is further discussed in Section~\ref{sec:tightness}. Note that for model $\mathrm{HYP}$-$\mathrm{SOC}$, depending on the optimal power consumption $p^*$, the approximation error and relaxation gap may offset each other.

A solution point $\mathrm{H4}$ with a zero relaxation gap can be recovered from point $\mathrm{H3}$ by reducing the power consumption of the electrolyzer from $p^*$ to $p^{**}$. Point $\mathrm{H4}$ is the only exact solution that has the same hydrogen production as point $\mathrm{H3}$.

In the following, we present the mathematical formulations for all the models.
\subsection{$\mathrm{HYP}$-$\mathrm{MIL}$}
\label{ssec: hyp_mil}
The state-of-the-art approximation of the hydrogen production curve \citep{powertech2023, YiZeng_PWL, MAROCCO_milp} follows a piecewise linearization approach as
\begingroup
\allowdisplaybreaks
\begin{subequations} 
\begin{align}
    & h_t = \sum_{s \in \mathcal{S}} (B_s \widehat{p}_{s,t} + C_s  \widehat{z}_{s,t}), & \forall\, & t \in \mathcal{T}, \label{eq:DA_lin_hydrogen_prod_MILP}\\
    & \underline{P}_s  \widehat{z}_{s,t} \leq \widehat{p}_{s,t} \leq  \overline{P}_s  \widehat{z}_{s,t}, & \forall\, & t \in \mathcal{T}, & \forall\, & s \in \mathcal{S},\label{eq:DA_p_seg_MILP}\\
    & z^{\mathrm{on}}_{t} = \sum_{s \in \mathcal{S}} \widehat{z}_{s,t}, & \forall\, & t \in \mathcal{T} ,\label{eq:DA_def_segments_MILP} \\
    & p_t = P^{\mathrm{sb}} z^{\mathrm{sb}}_{t} + \sum_{s \in \mathcal{S}} \widehat{p}_{s,t} & \forall\, & t \in \mathcal{T}, \label{eq:DA_power_consumption_tot_MILP}  \\
    & \widehat{z}_{s,t} \in \{0,1\}, & \forall\, & t \in \mathcal{T}, & \forall\, & s \in \mathcal{S}.\label{eq:DA_binaries_MILP}
\end{align}
\end{subequations} 
The hydrogen production on each linear segment $s$ is determined by \eqref{eq:DA_lin_hydrogen_prod_MILP}, where $B_s$ and $C_s$ denote the slope and intercept of the underlying segment, respectively. These coefficients have to be determined ex-ante by choosing fixed linearization points on the original nonlinear hydrogen production curve and performing a linear interpolation between them (see Figure~\ref{fig:3models}). In addition, variable $\widehat{p}_{s,t}$ is the power consumption corresponding to segment $s$ in time step $t$. Binary variable $\widehat{z}_{s,t}$ indicates whether the electrolyzer operates on segment $s$ in time step $t$. Constraint~\eqref{eq:DA_p_seg_MILP} enforces that the power consumption on each segment is between the lower and upper bounds $\underline{P}_s$ and $\overline{P}_s$, respectively, if the segment is active. Constraint~\eqref{eq:DA_def_segments_MILP} ensures that only one segment can be active when the electrolyzer is in on state. The power consumption of the electrolyzer is then defined by \eqref{eq:DA_power_consumption_tot_MILP} depending on the operational state. For $|\mathcal{S}| = 1$, model $\mathrm{HYP}$-$\mathrm{MIL}$ represents a linear hydrogen production curve, which does not require binary variables, as done in \cite{VARELA20219303}.
The interested reader is referred to \cite{powertech2023} for a more detailed explanation of the approximation of the hydrogen production curve based on piecewise linearization and a discussion on the impact of the number of segments.

\subsection{$\mathrm{HYP}$-$\mathrm{L}$}
The piecewise linearization of the hydrogen production curve \eqref{eq:DA_lin_hydrogen_prod_MILP}--\eqref{eq:DA_binaries_MILP} becomes computationally challenging with an increasing number of segments due to the required binary variables. A natural idea would then be to relax \eqref{eq:DA_lin_hydrogen_prod_MILP}, which allows the removal of associated binary variables. By this, \eqref{eq:DA_lin_hydrogen_prod_MILP}--\eqref{eq:DA_binaries_MILP} reduce to
\begin{subequations} 
\begin{align}
    & h_t \leq B_s \widetilde{p}_{t}+ C_s z^{\mathrm{on}}_{t}, & \forall\, & t \in \mathcal{T}, & \forall\, & s \in \mathcal{S},  \label{eq:DA_lin_hydrogen_prod_LP}\\
    & P^{\mathrm{min}}z^{\mathrm{on}}_{t} \leq \widetilde{p}_{t} \leq  P^{\mathrm{max}}z^{\mathrm{on}}_{t}, & \forall\, & t \in \mathcal{T},\label{eq:DA_p_hat_LP}\\
    & p_t = \widetilde{p}_{t} + P^{\mathrm{sb}} z^{\mathrm{sb}}_{t}, & \forall\, & t \in \mathcal{T}. \label{eq:DA_power_consumption_tot_LP}
\end{align}
\end{subequations}
The power consumption of the electrolyzer in on state is now given by $\widetilde{p}_t$. Constraint~\eqref{eq:DA_lin_hydrogen_prod_LP} is an intersection of hypographs of concave functions and are therefore convex~\citep{Boyd2004}. The set of constraints~\eqref{eq:DA_lin_hydrogen_prod_LP}--\eqref{eq:DA_power_consumption_tot_LP} is equivalent to \eqref{eq:DA_lin_hydrogen_prod_MILP}--\eqref{eq:DA_power_consumption_tot_MILP} when~\eqref{eq:DA_lin_hydrogen_prod_LP} is binding at the optimal solution. 
In contrast to piecewise linearization, this formulation is computationally efficient even for a comparatively high number of segments.

\subsection{$\mathrm{HYP}$-$\mathrm{SOC}$}
By looking into the experimental hydrogen production curve in Figure~\ref{fig:curves}, one may hypothesize that it has a quadratic shape. This is further supported by a semi-linear power-to-hydrogen conversion efficiency for power consumption levels higher than $P^{\mathrm{\eta, max}}$. Accordingly, we approximate the experimental hydrogen production curve $\mathrm{HYP}$-$\mathrm{X}$ by a second-order polynomial:
\begin{align}
    & h_{t} = A p_{t}^2 + B p_{t} + C, & \forall\, & t \in \mathcal{T}, \label{eq:DA_hydrogen_nonconvex}
\end{align}
where $A < 0$, $B > 0 $, and $C < 0$. As mentioned earlier, the second-order polynomial can be straightforwardly fitted to operational data of the electrolyzer, if available. This is a great advantage over models $\mathrm{HYP}$-$\mathrm{MIL}$ and $\mathrm{HYP}$-$\mathrm{L}$, as it does not require choosing the number and the location of linearization points.

Constraint~\eqref{eq:DA_hydrogen_nonconvex} is a nonconvex quadratic equality constraint. With binary variables required for modeling the operational states of the electrolyzer, the resulting optimization model would be a mixed-integer nonlinear programming (MINLP) problem, which is generally hard to solve with existing off-the-shelf solvers, even to locally optimal solutions.
In contrast to existing approaches in the literature, e.g., as in \cite{VARELA20219303}, we propose using a relaxed version of~\eqref{eq:DA_hydrogen_nonconvex}, which can be amended to include the operational states of the electrolyzer as
\begin{subequations} 
\begin{align}
    & h_{t} \leq A \widetilde{p}_{t}^2 + B \widetilde{p}_{t} + C z^{\mathrm{on}}_{t}, & \forall\, & t \in \mathcal{T}, \label{eq:DA_hydrogen_prod_SOC} \\
    & P^{\mathrm{min}}z^{\mathrm{on}}_{t} \leq \widetilde{p}_{t} \leq  P^{\mathrm{max}}z^{\mathrm{on}}_{t}, & \forall\, & t \in \mathcal{T},\label{eq:DA_p_hat_SOC}\\
    & p_t = \widetilde{p}_{t} + P^{\mathrm{sb}} z^{\mathrm{sb}}_{t}, & \forall\, & t \in \mathcal{T}. \label{eq:DA_power_consumption_tot_SOC}  
\end{align}
\end{subequations} 
For every time step $t$ that the electrolyzer is on, i.e., $z_{t}^{\mathrm{on}} = 1$, \eqref{eq:DA_hydrogen_prod_SOC} enforces a convex quadratic inequality constraint for $A < 0$. This constraint can be reformulated into an efficiently solvable rotated second-order cone (SOC) constraint, as shown in Section~OC.2 of the online companion \citep{Online_Companion}. In the following, we stick to the convex quadratic form as we find it to be more intuitive here.
The set of constraints~\eqref{eq:DA_hydrogen_prod_SOC}--\eqref{eq:DA_power_consumption_tot_SOC} is equivalent to \eqref{eq:DA_hydrogen_nonconvex} when~\eqref{eq:DA_hydrogen_prod_SOC} is binding for every time step $t$ at the optimal solution. We then say that the relaxation is exact.
If \eqref{eq:DA_hydrogen_prod_SOC} is nonbinding, i.e., there is a difference between the right and left-hand side of constraint~\eqref{eq:DA_hydrogen_prod_SOC}, then the relaxation gap is nonzero. This is illustrated by the discrepancy between points $\mathrm{H2}$ and $\mathrm{H3}$ in Figure~\ref{fig:3models}.

Accounting for binary variables needed to model the operational states of the electrolyzer, the resulting problem is a mixed-integer second-order cone programming (MISOCP) problem, which is efficiently solvable by existing off-the-shelf solvers like Gurobi, Mosek, and CPLEX, that directly support convex quadratic constraints.

\subsection{Summary}
An overview of the three models is given in Table~\ref{tab:model_overview}. Hereafter, we use terms $\mathrm{HYP}$-$\mathrm{MIL}$, $\mathrm{HYP}$-$\mathrm{L}$, and $\mathrm{HYP}$-$\mathrm{SOC}$, not only to refer to the hydrogen production curve but also to the resulting scheduling problem of the hybrid power plant.

The objective function \eqref{eq:DA_objective} is common to all three models. 
We group the constraints into three groups. The first group consists of all constraints related to the hybrid power plant as a system, which are all linear. The second group includes all constraints related to the operational states of the electrolyzer, requiring binary variables. The first and second groups are common to all three models. The third group depends on the model of the hydrogen production curve that replaces constraint~\eqref{eq:DA_hydrogen_general}, i.e., the state-of-the-art piecewise linear approximation ($\mathrm{HYP}$-$\mathrm{MIL}$), a linear relaxation counterpart ($\mathrm{HYP}$-$\mathrm{L}$), or the proposed conic relaxation ($\mathrm{HYP}$-$\mathrm{SOC}$).
Due to the binary variables $\mathbf{z}$ for modeling the operational states, the resulting scheduling problem of the hybrid power plant becomes either a mixed-integer linear programming (MILP) or a MISOCP problem.
When neglecting the operational states, models $\mathrm{HYP}$-$\mathrm{L}$ and $\mathrm{HYP}$-$\mathrm{SOC}$ reduce to a linear programming (LP) or a second-order cone programming (SOCP) problem, respectively. This is not the case for model $\mathrm{HYP}$-$\mathrm{MIL}$, which requires additional binary variables for the piecewise linearization of the hydrogen production curve.

We define the set of variables related to the hydrogen production curve $\mathbf{y}$ as $\mathbf{y}^{\mathrm{L}}=\mathbf{y}^{\mathrm{SOC}}=\{\widetilde{p}_t\}$, and $\mathbf{y}^{\mathrm{MIL}}=\{\widehat{p}_{s,t}, \widehat{z}_{s,t}\}$ for $\mathrm{HYP}$-$\mathrm{L}$, $\mathrm{HYP}$-$\mathrm{SOC}$, and $\mathrm{HYP}$-$\mathrm{MIL}$, respectively.

\begin{table}
\centering
\caption{Summary of constraints for the three different models. The objective function \eqref{eq:DA_objective} is common to all models.}
\vspace{0.1cm}
\label{tab:model_overview}
\begin{tabular}{l|ccc}
& \textbf{$\mathrm{HYP}$-$\mathrm{MIL}$} & \textbf{$\mathrm{HYP}$-$\mathrm{L}$} & \textbf{$\mathrm{HYP}$-$\mathrm{SOC}$} \\
\hline
\begin{tabular}[c]{@{}l@{}} Hybrid \\ power \\ plant \end{tabular} & 
\multicolumn{3}{c}{
\begin{tabular}[c]{@{}l@{}}
Power balance: \eqref{eq:DA_balance} \\
Grid exchange: \eqref{eq:DA_green_hydrogen} \\ 
Hydrogen demand: \eqref{eq:DA_hydrogen_demand}-- \eqref{eq:DA_non_negative_hydrogen}
\end{tabular}} \\ \hline
\begin{tabular}[c]{@{}l@{}} Electrolyzer \\ states \end{tabular} & 
\multicolumn{3}{c}{
\begin{tabular}[c]{@{}l@{}} 
Mutual exclusiveness: \eqref{eq:DA_mutual_states} \\
Power consumption limits: \eqref{eq:DA_consumption_max}--\eqref{eq:DA_consumption_min}  \\
Startup and state binaries: \eqref{eq:DA_start_up}--\eqref{eq:DA_start_up_1},  \eqref{eq:DA_state_binaries}\\
\end{tabular}} \\
\hline
\begin{tabular}[c]{@{}l@{}} Hydrogen \\ production \end{tabular} &
\eqref{eq:DA_lin_hydrogen_prod_MILP}--\eqref{eq:DA_binaries_MILP} &
\eqref{eq:DA_lin_hydrogen_prod_LP}--\eqref{eq:DA_power_consumption_tot_LP} &
\eqref{eq:DA_hydrogen_prod_SOC}--\eqref{eq:DA_power_consumption_tot_SOC}
 \\ \hline
\begin{tabular}[c]{@{}l@{}} Model \\ type \end{tabular}  &
MILP & MILP & MISOCP \\
\hline
\end{tabular}
\end{table}

%% file: 03A_Tightness.tex
\section{On the Exactness of the Conic Relaxation}
\label{sec:tightness}
Models $\rm{HYP}$-$\rm{L}$ and $\rm{HYP}$-$\rm{SOC}$ are exact when inequality constraints \eqref{eq:DA_lin_hydrogen_prod_LP} and \eqref{eq:DA_hydrogen_prod_SOC}, respectively, are binding at the optimal solution.
In the following, we derive sufficient conditions for model $\rm{HYP}$-$\rm{SOC}$ to be exact as well as necessary and sufficient conditions for it to be inexact. Similar analytical results can be obtained for $\rm{HYP}$-$\rm{L}$. Focusing on $\rm{HYP}$-$\rm{SOC}$, if \eqref{eq:DA_hydrogen_prod_SOC} is nonbinding at optimum, i.e., the relaxation gap is nonzero, an intuitive interpretation is that a fraction of hydrogen produced is being wasted (i.e., the difference between points $\mathrm{H2}$ and $\mathrm{H3}$ in Figure~\ref{fig:3models}) or that an unnecessarily high amount of power is being consumed (i.e., the difference between points $\mathrm{H3}$ and $\mathrm{H4}$ in Figure~\ref{fig:3models}).

In the following, we always assume that a solution to problem $\rm{HYP}$-$\rm{SOC}$ exists, as $p=0, f=W$, and $z^{\mathrm{off}}=1$ is trivial. Recall we assume that the hydrogen price is positive, such that the objective function is strictly increasing in hydrogen production, as stated in Section~\ref{sec:prelim}. Now, we provide analytical results for the exactness of relaxation~\eqref{eq:DA_hydrogen_prod_SOC}.

\begin{theorem}\label{Theorem_1}
If the maximum total hydrogen production constraint~\eqref{eq:DA_hydrogen_demand} is nonbinding at optimum, then the relaxation~\eqref{eq:DA_hydrogen_prod_SOC} is exact.
\end{theorem}

\textit{Proof:} See \ref{sec:proof_theorem_1}. 

Intuitively, if the hydrogen price is positive, there is no incentive to waste hydrogen as long as it can be sold to the demand. The hydrogen production corresponds to point $\mathrm{H2}$ in Figure~\ref{fig:3models}. This theorem applies to a wide range of relevant use cases wherein the total hydrogen production is unconstrained. This could be the case wherein hydrogen production facilities are located next to large hydrogen consumers or a direct connection to hydrogen pipeline infrastructure or large storage facilities exists. In some real-life applications, however, the maximum total hydrogen production is constrained, e.g., by the capacity of available tube trailers for hydrogen transportation. For that case, we derive another theorem providing a sufficient condition for exactness based on electricity prices being positive.
\begin{theorem}\label{Theorem_2}
Suppose the maximum total hydrogen production~\eqref{eq:DA_hydrogen_demand} is binding at optimum for sub-period $n$. If the power prices are positive $\lambda_t > 0, \forall t \in \mathcal{H}_n$, then the relaxation~\eqref{eq:DA_hydrogen_prod_SOC} is exact.
\end{theorem}

\textit{Proof:} See \ref{sec:proof_theorem_2}. 

Intuitively, when the power price is positive, Theorem~\ref{Theorem_2} implies that increasing the power consumption of the electrolyzer is unprofitable if additional hydrogen cannot be sold. To complete our analyses, we now focus on cases with nonpositive electricity prices, which constitute a highly profitable business case for an electrolyzer.
When buying power from the grid and spilling wind is not allowed, there is a monetary incentive for the electrolyzer to maximize its consumption of local wind power production, even if the hydrogen demand is already satisfied.

\begin{theorem}\label{Theorem_3}
Let $\mathcal{T}_{n}^- \subseteq \mathcal{H}_{n}$ be the subset of hours in sub-period $n$, such that $\lambda_t \leq 0, \forall t \in \mathcal{T}_{n}^-$. If the total hydrogen production in those hours equals the maximum demand, i.e., $\sum_{t \in \mathcal{T}_{n}^-} h_t = D_{n}^{\rm{max}}$, such that \eqref{eq:DA_hydrogen_demand} is binding at optimum for sub-period ${n}$, then there exists at least one time step $t \in \mathcal{T}_{n}^-$ for which the relaxation~\eqref{eq:DA_hydrogen_prod_SOC} is inexact.
\end{theorem}

\textit{Proof:} See \ref{sec:proof_theorem_3}. 

Intuitively, when power prices are negative, it is most profitable to increase the power consumption of the electrolyzer as much as possible. If the corresponding hydrogen production exceeds the maximum total demand level, it is wasted without additional cost. The wasted hydrogen production equals the difference between points $\mathrm{H2}$ and $\mathrm{H3}$ in Figure~\ref{fig:3models}. Due to the opportunity cost related to the negative power prices, this is more profitable than reducing the power consumption to $\mathrm{H4}$.
While Theorem~\ref{Theorem_3} states that relaxation~\eqref{eq:DA_hydrogen_prod_SOC} is inexact when it is economically most profitable for the hybrid power plant, the necessary conditions are hardly met in practice. We demonstrate this based on a realistic case study in Section~\ref{sec:results}.
Note that the theorems stated here, including Theorem~\ref{Theorem_3}, extend to the case where the electrolyzer is allowed to buy electricity from the grid and wind spillage is allowed.

The necessary condition for inexactness of relaxation~\eqref{eq:DA_hydrogen_prod_SOC} stated in Theorem~\ref{Theorem_3} can be checked a priori. For each sub-period $n$, one can evaluate in advance if the total maximum possible hydrogen production during hours with nonpositive prices is greater than or equal to the maximum demand, i.e.,
\begin{align}
    & \sum_{t \in \mathcal{T}_n^-} A \overline{P}_{t}^2 + B \overline{P}_{t} + C \geq D^{\mathrm{max}}_{n}, & \forall\, & n \in \mathcal{N}, \label{eq:cc}
\end{align}
where $\overline{P}_{t} = \min\{W_t, P^{\mathrm{max}}\}, \ \forall t \in \mathcal{T}$. It follows from Theorem~\eqref{Theorem_3} that if condition~\eqref{eq:cc} is not fulfilled, then the relaxation is exact. The converse is not necessarily true.
If the relaxation~\eqref{eq:DA_hydrogen_prod_SOC} is inexact, model $\rm{HYP}$-$\rm{SOC}$ may obtain a solution like $\mathrm{H3}$ in Figure~\ref{fig:3models}. This might include a sub-optimal solution for the binary variables related to the operational states of the electrolyzer.

\begin{figure}[t]
\centering
\includegraphics[width=0.8\linewidth]{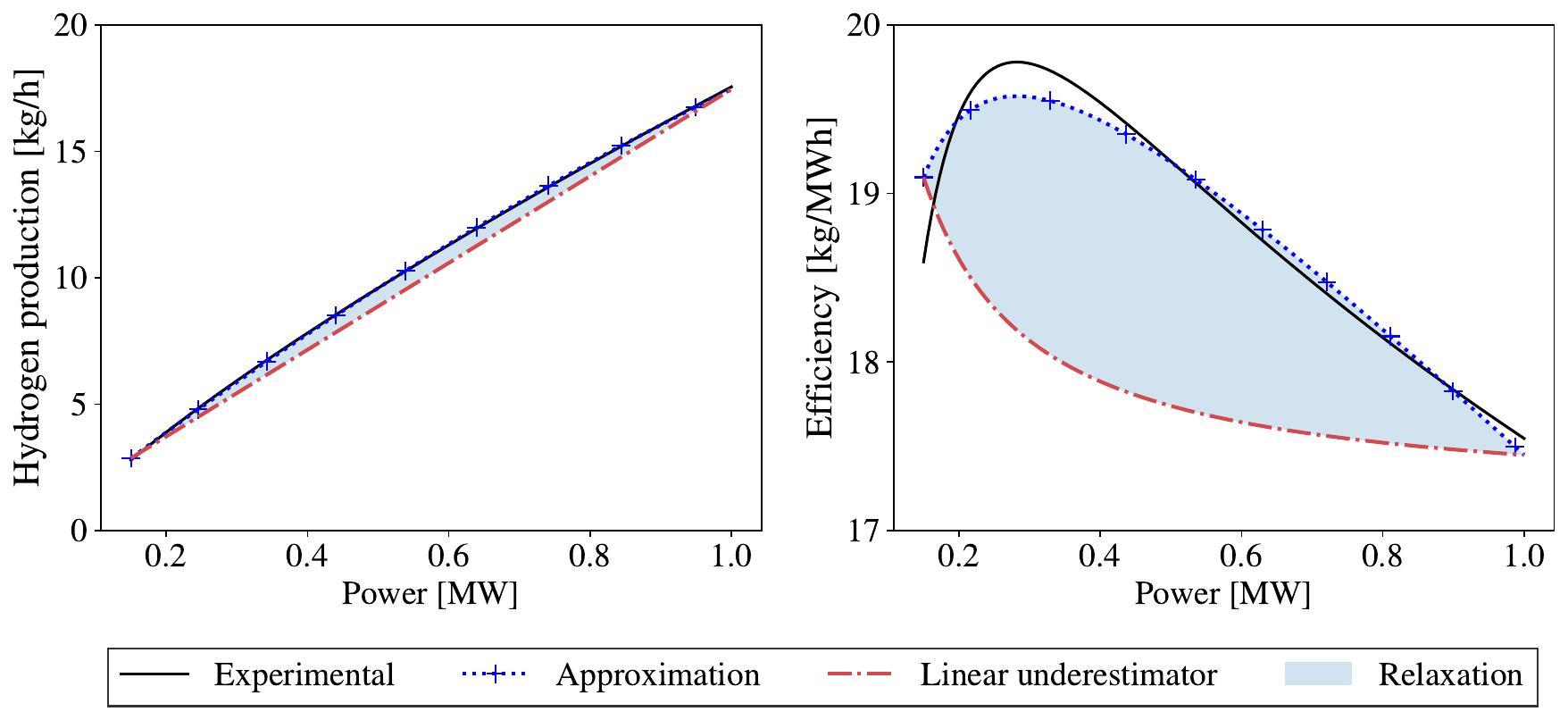}
\caption{Intersection of the feasible regions defined by the hydrogen production curve $\rm{HYP}$-$\rm{SOC}$ and the corresponding linear underestimator for the hydrogen production (left) and efficiency (right).} 
\label{fig:linear_underestimator}
\end{figure}

In cases where it cannot be guaranteed a priori by using the heuristic in \eqref{eq:cc} that the relaxation gap will be zero, the feasible region should be tightened. One option is to enforce a linear underestimator \citep{Taylor2021} to the chosen model of the hydrogen production curve. For model $\rm{HYP}$-$\rm{SOC}$, this can be done by adding the following constraint:
\begin{align}
    & h_t \geq Q_1 \widetilde{p}_t + Q_0, & \forall\, & t \in \mathcal{T}, \label{eq:DA_linear_underestimator_SOC}
\end{align}
where coefficients $Q_0$ and $Q_1$ can be derived from a linear interpolation of the points corresponding to the minimum and maximum power consumption. The linear underestimator~\eqref{eq:DA_linear_underestimator_SOC} substantially reduces the feasible region of model $\rm{HYP}$-$\rm{SOC}$, as shown in Figure~\eqref{fig:linear_underestimator}.
In the case considered here, the theoretical magnitude of the relaxation gap for a single time step is reduced from around 17.55 kg/h to 0.7 kg/h.

A nonzero relaxation gap, e.g., point $\mathrm{H3}$ in Figure~\eqref{fig:3models}, implies that the solution will be physically infeasible for the electrolyzer system in practice. In that case, a physically feasible solution should be restored a posteriori. For the hybrid power plant considered here, when using the linear underestimator~\eqref{eq:DA_linear_underestimator_SOC}, this can always be done by reducing the power consumption of the electrolyzer from $p^{*}$ at point $\mathrm{H3}$ to $p^{**}$ at point $\mathrm{H4}$, as shown in Figure~\eqref{fig:3models}. This might lead to sub-optimal solutions in more complex setups, e.g., when considering additional chemical downstream processes. In those cases, other techniques for recovering a feasible solution can be applied, e.g., the convex-concave procedure as done by \citet{Taylor2022} for biochemical processes or \citet{Wei2017} for AC optimal power flow. Other options are discussed by \citet{Venzke2020} in the context of AC optimal power flow.

%% file: 04_Results.tex
\section{Numerical Study}
\label{sec:results}
We consider a hybrid power plant consisting of a $1$-MW electrolyzer and a $2$-MW wind farm. Hourly wind capacity factors for year $2019$ are obtained from the Renewable.ninja web platform \citep{renewables.ninja_2} for a wind farm located in Eastern Denmark. The day-ahead electricity prices in the same year for the corresponding bidding area ($\mathrm{DK2}$) are taken from \cite{entso-e}. The electrolyzer has a minimum operating power of $P^{\mathrm{min}}$ = $0.15$ MW and a standby power consumption of $P^{\mathrm{sb}}$ = $0.01$ MW. We assume a cold-startup cost of $K^{\mathrm{su}}$ = \euro$50$, borrowed from \cite{VARELA20219303}. Hydrogen is sold at a fixed price of $\chi$ = \euro$2.1$/kg. Recall that electricity must not be purchased from the grid.

The public repository \citep{Online_Companion} contains the input data and code implementation in the Pyomo package \citep{hart2011pyomo,bynum2021pyomo} for Python, where optimization models have been solved by the Gurobi solver \citep{gurobi}.

\subsection{Exactness of the Proposed Conic Relaxation}
\label{sec:case_study_tightness}
\begin{figure*}[h]
\centering
\includegraphics[width=1\linewidth]{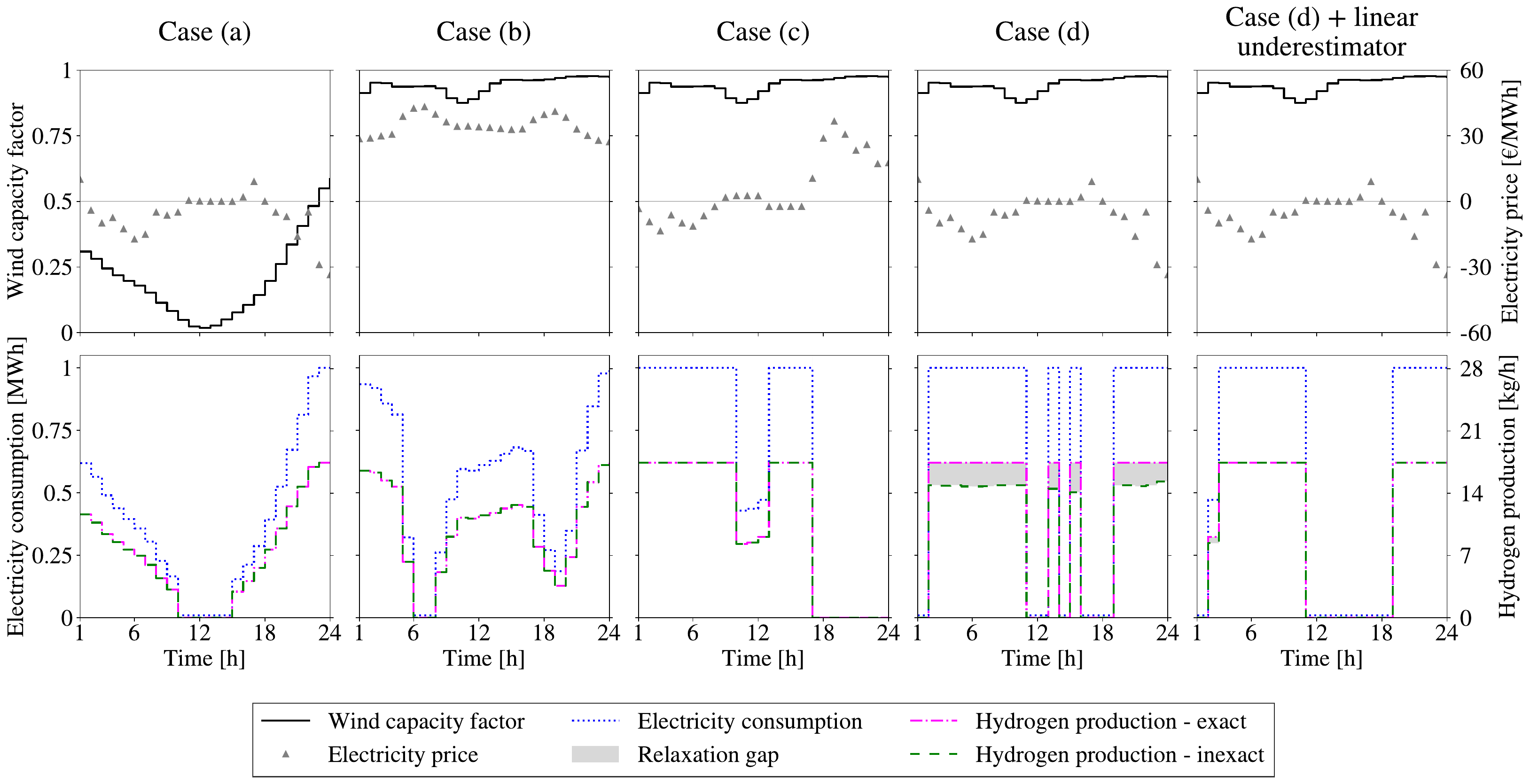}%
\caption{Top row: Wind capacity factor (left axis) and day-ahead power price (right axis) profiles for four illustrative days. Bottom row: Optimal electrolyzer power consumption (blue) and hydrogen production (green) for model $\mathrm{HYP}$-$\mathrm{SOC}$. The pink line shows the hydrogen production corresponding to a zero relaxation gap, which is illustrated by the shaded area (grey). When the relaxation~\eqref{eq:DA_hydrogen_prod_SOC} is exact, i.e., the gap is zero, the green and pink lines coincide. The last column shows the result for Case~(d) when adding the linear underestimator \eqref{eq:DA_linear_underestimator_SOC} to model $\mathrm{HYP}$-$\mathrm{SOC}$.} 
\label{fig:tightness}
\end{figure*}

As discussed in Section~\ref{sec:tightness}, the exactness of the relaxation in $\mathrm{HYP}$-$\mathrm{SOC}$ depends on the power prices and total hydrogen production. The latter is affected by the wind power availability. In the remainder of the paper, we use the term operating conditions to indicate a combination of power prices, wind availability, and type of hydrogen demand (i.e., constrained or unconstrained). To validate our analytical results in Section~\ref{sec:tightness}, we consider four cases with different profiles of wind power availability and power prices:
\begin{enumerate}[label = {Case~(\alph*):},align=left]
    \item Mostly negative hourly power prices and low wind power availability,
    \item Solely positive hourly power prices and comparatively high wind power availability compared to that in Case~(a),
    \item Some negative hourly power prices and comparatively high wind power availability as in Case~(b),
    \item Mostly negative hourly power prices as in Case~(a) and comparatively high wind power availability as in Case~(b).
\end{enumerate}
For each of the four cases, we build an illustrative case study for a time period of one day ($24$ hours). For that, we select different combinations of two realistic wind profiles (low and high wind) and three power price profiles (solely positive, some negative, and mainly negative prices).
The data is visualized in the top row of Figure~\ref{fig:tightness}. Among the four illustrative days, the one representing Case~(c) is the only one that combines wind and price profiles of the same day in $2019$. Note that the two price profiles in Cases~(c) and (d) correspond to the two days in $2019$ with the highest number of hours with negative power prices. In total, the dataset for $2019$ includes $95$ hours with negative electricity prices, spread across 20 days. Days that belong to Case~(d) were only observed twice in $2019$.

To create a situation where the maximum daily hydrogen production constraint~\eqref{eq:DA_hydrogen_demand} becomes binding, we choose a rather low hydrogen demand of $D^{\mathrm{max}} = 252.7\,\mathrm{kg}$ for all cases. This value corresponds to the amount of hydrogen that is produced when the electrolyzer is operated in full-load operation $60$\% of the time.
The optimal power consumption, hydrogen production, and the relaxation gap for Cases (a)-(d) are shown in the bottom row of Figure~\ref{fig:tightness}. Recall that by relaxation gap, we refer to the discrepancy between the right and left-hand sides of constraint~\eqref{eq:DA_hydrogen_prod_SOC}, which is illustrated by points $\mathrm{H2}$ and $\mathrm{H3}$ in Figure~\ref{fig:3models}. Note that all results on the exactness of the relaxation of the hydrogen production curve presented in this subsection can be similarly obtained for model $\mathrm{HYP}$-$\mathrm{L}$.

In Case~(a), despite frequent negative power prices, there is not enough wind power available over the day such that the total hydrogen demand~\eqref{eq:DA_hydrogen_demand} reaches the upper limit.
The comparably higher wind power availability in Case~(b) in combination with low but positive power prices results in a binding maximum total hydrogen demand.
For both Cases (a) and (b), it follows immediately from Theorem~\ref{Theorem_1} and Theorem~\ref{Theorem_2}, respectively, that the relaxation $\mathrm{HYP}$-$\mathrm{SOC}$ is exact. 

The maximum total hydrogen production is binding in Case~(c) too. It is, however, non-binding for the 12 hours with negative prices only. Hence, it does not fulfill the necessary condition for the inexactness of relaxation~\eqref{eq:DA_hydrogen_prod_SOC} as stated in Theorem~\ref{Theorem_3}. Note that the operation during hours with negative prices is preferred even if operating in hours with positive prices is profitable. Case~(d) has an even higher number of hours with negative prices (17 hours), such that the hydrogen demand becomes binding during the operation in those hours only. This satisfies the necessary and sufficient conditions stated in Theorem~\ref{Theorem_3}. As a consequence, relaxation~\eqref{eq:DA_hydrogen_prod_SOC} is inexact at optimum, leading to a non-zero relaxation gap equal to the shaded area in Figure~\ref{fig:tightness}.\footnote{Note that for model $\mathrm{HYP}$-$\mathrm{SOC}$ the solver reports a solution where all time steps have a non-zero feasibility gap. This is not necessarily always true, as multiple optimal solutions exist in this case. When solving the same case study with the linear relaxation $\mathrm{HYP}$-$\mathrm{L}$, the aggregated relaxation gap is fully allocated to the minimum possible number of time steps.}

Recall that Case~(d) represents a rare combination of prolonged hours with negative power prices (17 hours), excessive wind availability ($95$\% average wind capacity factor), and restrictive maximum hydrogen demand ($60$\% of the total full-load production).
For all other combinations of those three factors, which we define as prevalent operating conditions, the relaxation is exact.
While negative power prices are rarely observed in electricity markets at the moment, they might occur more frequently in the future. Nonetheless, this will not compromise the applicability of the proposed relaxations in terms of exactness, as it additionally requires a constrained maximum production limit (Theorem~\ref{Theorem_3}). Recall that we assume a highly restrictive maximum demand limit set by tube trailers with fixed schedules to illustrate our mathematical findings on the exactness of the relaxations. For the future, we envision that electrolyzers will be placed close to pipeline infrastructure, large storage facilities, or industrial consumers, in which the maximum hydrogen demand can be considered unconstrained.

In the last column of Figure~\ref{fig:tightness}, Case~(d) is solved after adding the linear underestimator~\eqref{eq:DA_linear_underestimator_SOC} to model $\mathrm{HYP}$-$\mathrm{SOC}$.
Compared to Case~(d), this reduces the number of hours with an inexact relaxation from 17 to 1 and the total relaxation gap from 43 kg to 0.7 kg of hydrogen. It can further be noticed that the electrolyzer is now operated in standby state instead of on state in hours 13 and 15.
For $\mathrm{HYP}$-$\mathrm{SOC}$ without linear underestimator, this was possible because the hydrogen production in on state is bounded below by zero, as can be seen in Figure~\ref{fig:3models}. After adding the linear underestimator~\eqref{eq:DA_linear_underestimator_SOC}, the minimum hydrogen production in on state is 2.9 kg/h, as can be seen in Figure~\ref{fig:linear_underestimator}. Since the maximum hydrogen demand~\eqref{eq:DA_hydrogen_demand} is binding, the minimum hydrogen production cannot be accommodated when including the linear underestimator. Consequently, the electrolyzer must be turned into standby or off state. This indicates that an inexact relaxation may lead to a sub-optimal solution for the binary variables related to the operational states of the electrolyzer, compared to solving the original MINLP including nonconvex constraint~\eqref{eq:DA_hydrogen_nonconvex} directly.

\subsection{Comparison of the Solution Quality}
\label{sec:results_quality}
This section compares the models proposed in Section~\ref{sec:model} in terms of profit, dispatch decisions, and hydrogen production.
Based on the heuristic proposed in Equation~\eqref{eq:cc}, we found that for the year $2019$, the proposed relaxations are exact for a maximum hydrogen demand above $D^{\mathrm{max}} = 296.3\,\mathrm{kg}$, corresponding to the hydrogen produced when the electrolyzer runs in full-load operation $70.4$\% of the time. For the relaxations to be inexact for more than $2$ days, an unrealistically low maximum hydrogen demand below $25$\% is necessary. This emphasizes that the proposed relaxations are exact under prevalent operating conditions.
Therefore, we choose a maximum daily hydrogen demand corresponding to 90\%, which ensures that the relaxations $\mathrm{HYP}$-$\mathrm{L}$ and $\mathrm{HYP}$-$\mathrm{SOC}$ are always exact. This implies that $\mathrm{HYP}$-$\mathrm{L}$ results in the same dispatch decisions as $\mathrm{HYP}$-$\mathrm{MIL}$. Therefore, we do not explicitly report the results for model $\mathrm{HYP}$-$\mathrm{L}$. 

We use the term $\mathrm{HYP}$-$\mathrm{MIL}i$ for $i \in \{1,2,10,24\}$ to refer to model $\mathrm{HYP}$-$\mathrm{MIL}$ with $i$ linearization segments. Model $\mathrm{HYP}$-$\mathrm{MIL}24$, i.e., model $\mathrm{HYP}$-$\mathrm{MIL}$ with 24 segments, is used as a benchmark in the following.
A detailed explanation of how we choose the location of linearization points is given in Section~OC.3 of the online companion \citep{Online_Companion}. The interested reader is further referred to \cite{MAROCCO_milp} for a different choice of linearization points.

\begin{figure}[t]
\centering
\includegraphics[width=0.6\linewidth]{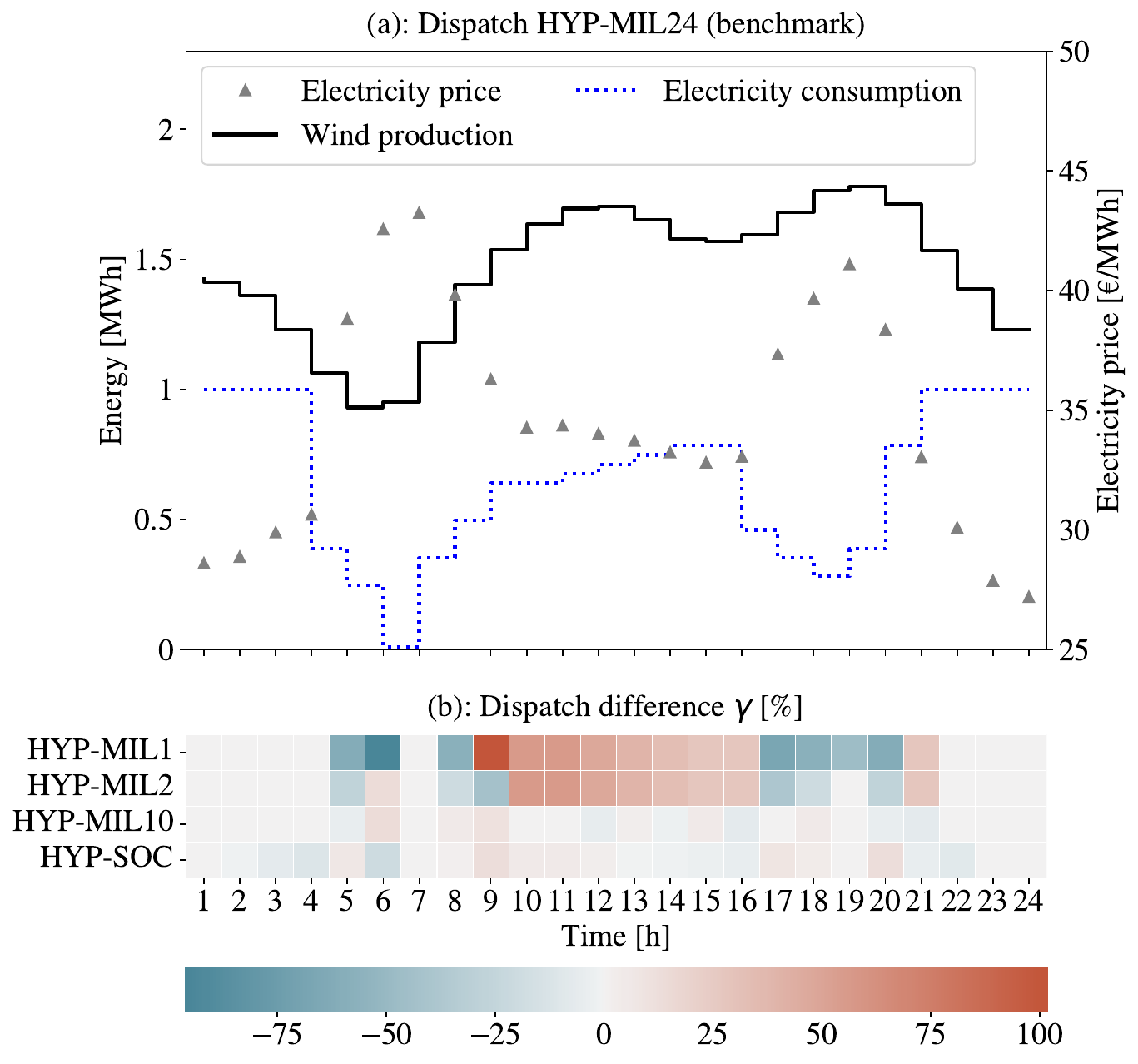}
\caption{(a): Optimal electrolyzer power consumption for the $\mathrm{HYP}$-$\mathrm{MIL}24$ benchmark on an illustrative day. (b): Hourly relative difference $\gamma_t$ of the optimal electrolyzer power consumption compared to the benchmark.} 
\label{fig:1day_comparison}
\end{figure}

We first compare the different models in terms of dispatch decisions for a selected day in $2019$, where the hydrogen production is not limited by wind power availability.
Figure~\ref{fig:1day_comparison}(a) shows the wind profile, electricity prices, and optimal electrolyzer power consumption of $\mathrm{HYP}$-$\mathrm{MIL}24$, which is chosen as a benchmark. For each model formulation, the hourly relative difference of the optimal power consumption compared to the benchmark, denoted as $\gamma_t$, is defined as
\begin{align}
    & \gamma_t = \frac{p_t^*-p_t^\mathrm{*,MIL24}}{p_t^\mathrm{*,MIL24}}, & \forall\, & t \in \mathcal{T},
\end{align}
where $p_t^*$ is the optimal power consumption of the electrolyzer in hour $t$ obtained from the underlying model and $p_t^\mathrm{*,MIL24}$ is that of the benchmark. The relative difference $\gamma_t$ is shown for all models in Figure~\ref{fig:1day_comparison}(b). For models $\mathrm{HYP}$-$\mathrm{MIL1}$ and $\mathrm{HYP}$-$\mathrm{MIL2}$, the relative difference is comparatively large in hours where the electrolyzer operates at partial loading. This is the case when the electricity price is in the range of $31$ to $43$ \euro/MWh (cf. Appendix in \cite{powertech2023}). The average daily difference, i.e., $\bar{\gamma} = (\sum_t|\gamma_t|)/24$ for $\mathrm{HYP}$-$\mathrm{MIL1}$ and $\mathrm{HYP}$-$\mathrm{MIL2}$ is $36$\% and $21$\%, respectively. In comparison, $\mathrm{HYP}$-$\mathrm{SOC}$ exhibits significantly better performance, with an average error of around $5$\%. 
To achieve a similar solution quality with $\mathrm{HYP}$-$\mathrm{MIL}$ (and $\mathrm{HYP}$-$\mathrm{L}$), at least $10$ segments are needed.\footnote{Note that a different choice of linearization points, e.g., based on the method presented in \cite{MAROCCO_milp}, may slightly affect the number of segments needed to achieve a similar solution accuracy between $\mathrm{HYP}$-$\mathrm{SOC}$ and $\mathrm{HYP}$-$\mathrm{MIL}$/$\mathrm{HYP}$-$\mathrm{L}$.}

To validate the robustness of our findings presented in Figure~\ref{fig:1day_comparison}, we solve the different models for the entire year $2019$ and compare the profit and dispatch results. To fairly compare the different models in terms of hydrogen production, an ex-post analysis is performed to take the model-specific approximation errors into account. The ex-post analysis consists of solving the underlying optimization problem, fixing the optimal power consumption of the electrolyzer, and then calculating the corresponding hydrogen production based on $\mathrm{HYP}$-$\mathrm{X}$, indicated by point $\mathrm{H1}$ in Figure~\ref{fig:3models}.
The reason for the ex-post analysis is that the model determines the choice of dispatch decisions (i.e., how much power should be sold to the grid and how much power should be consumed by the electrolyzer), while the actual hydrogen production depends on the electrolyzer physics and not on the approximated model that is adopted in the scheduling problem. The difference in the hydrogen production quantities and profits derived from the scheduling model and those from the ex-post analysis, therefore, serve as an indicator of the quality of the approximation of the experimental hydrogen production curve $\mathrm{HYP}$-$\mathrm{X}$. Recall that we cannot solve the scheduling model including $\mathrm{HYP}$-$\mathrm{X}$ directly, as it is not possible to derive a corresponding analytical function.
If it was possible to solve scheduling model $\mathrm{HYP}$-$\mathrm{X}$ directly, then the results of its ex-post analysis would be exactly the same as those derived from the optimization model. To replicate this idealized benchmark as best as possible, we choose the ex-post results of model $\mathrm{HYP}$-$\mathrm{MIL24}$ as a benchmark for our case studies.
Note that we do not account for the hydrogen production corresponding to the approximation error in the scheduling problem, which could potentially make the maximum hydrogen production constraint~\eqref{eq:DA_hydrogen_demand} binding.

Table~\ref{tab:1year_comparison} presents the annual profit, electricity sales, and hydrogen production, as a difference compared to the $\mathrm{HYP}$-$\mathrm{MIL24}$ benchmark. While the total annual profit is similar for all models, the share of hydrogen production and power sales for $\mathrm{HYP}$-$\mathrm{MIL1}$ and $\mathrm{HYP}$-$\mathrm{MIL2}$ differ substantially compared to the $\mathrm{HYP}$-$\mathrm{MIL24}$ benchmark. Due to the inaccurate approximation of the hydrogen production curve, the number of hours when it is profitable to produce hydrogen is considerably reduced, leading to approximately $14$\% and $7$\% less hydrogen production, respectively. The corresponding monetary loss is partially compensated by an increase in wind energy sold to the grid. Model $\mathrm{HYP}$-$\mathrm{SOC}$ shows higher accuracy in decision-making, with a difference in hydrogen production lower than $1$\% compared to the benchmark. Similarly to the one-day example, at least $10$ segments are necessary to achieve the same results with $\mathrm{HYP}$-$\mathrm{MIL}$ as with the conic model.

\begin{table}[b]
\centering
\caption{Difference in the annual profit, hydrogen production, and power sales compared to the $\mathrm{HYP}$-$\mathrm{MIL24}$ benchmark. Note that results for $\mathrm{HYP}$-$\mathrm{L}$ are not explicitly reported as they coincide with $\mathrm{HYP}$-$\mathrm{MIL}$ for the considered case study.}
\vspace{0.1cm}
\begin{tabular}{lccc}
\hline
                                 & Profit      & Hydrogen production & Power sales \\ \hline
$\mathrm{HYP}$-$\mathrm{MIL}$24  & -           & -                   & -           \\
$\mathrm{HYP}$-$\mathrm{MIL}$10  & -0.003\%    & -0.12\%             & 0.05\%    \\
$\mathrm{HYP}$-$\mathrm{MIL}$2   & -0.26\%     & -7.06\%             & 3.66\%    \\
$\mathrm{HYP}$-$\mathrm{MIL}$1   & -0.63\%     & -13.84\%            & 6.52\%    \\
$\mathrm{HYP}$-$\mathrm{SOC}$    & -0.01\%     & -0.89\%             & 0.48\%    \\\hline
\end{tabular}
\label{tab:1year_comparison}
\end{table}

\subsection{Comparison of the Computational Performance}
\label{sssec:computational}
We now compare the different models in terms of their solution time, and then analyze how it scales with the problem size. To do so, we extend the deterministic day-ahead scheduling problem proposed in Table~\ref{tab:model_overview} to a two-stage stochastic problem and run it for a single day using different numbers of scenarios $\omega \in \Omega$. We assume perfect foresight for the electricity prices but uncertainty in wind power production. In the second stage, the electrolyzer adjusts its power consumption to minimize the real-time imbalance cost based on the scenario-specific wind power realization. Owed to the fast dynamics of the electrolyzer, we further assume that it is able to change its operational states compared to the day-ahead schedule. The stochastic model formulation is reported in Section~OC.4 of the online companion \citep{Online_Companion}. We use the same price profile as in the case study reported in Figure~\ref{fig:1day_comparison}. For the wind power scenarios, the first $24$ hours of the dataset provided in \cite{Pinson} are used. The problem was solved using a High Performance Computing Cluster node with two AMD EPYC $7551$ processors clocking at $2$ GHz and using a maximum of $8$ threads. As this section focuses on the computational performance only, we do not further elaborate on the quality of the solution for the stochastic model and look solely at in-sample results.

\begin{figure}[]
\centering
\includegraphics[width=0.6\linewidth]{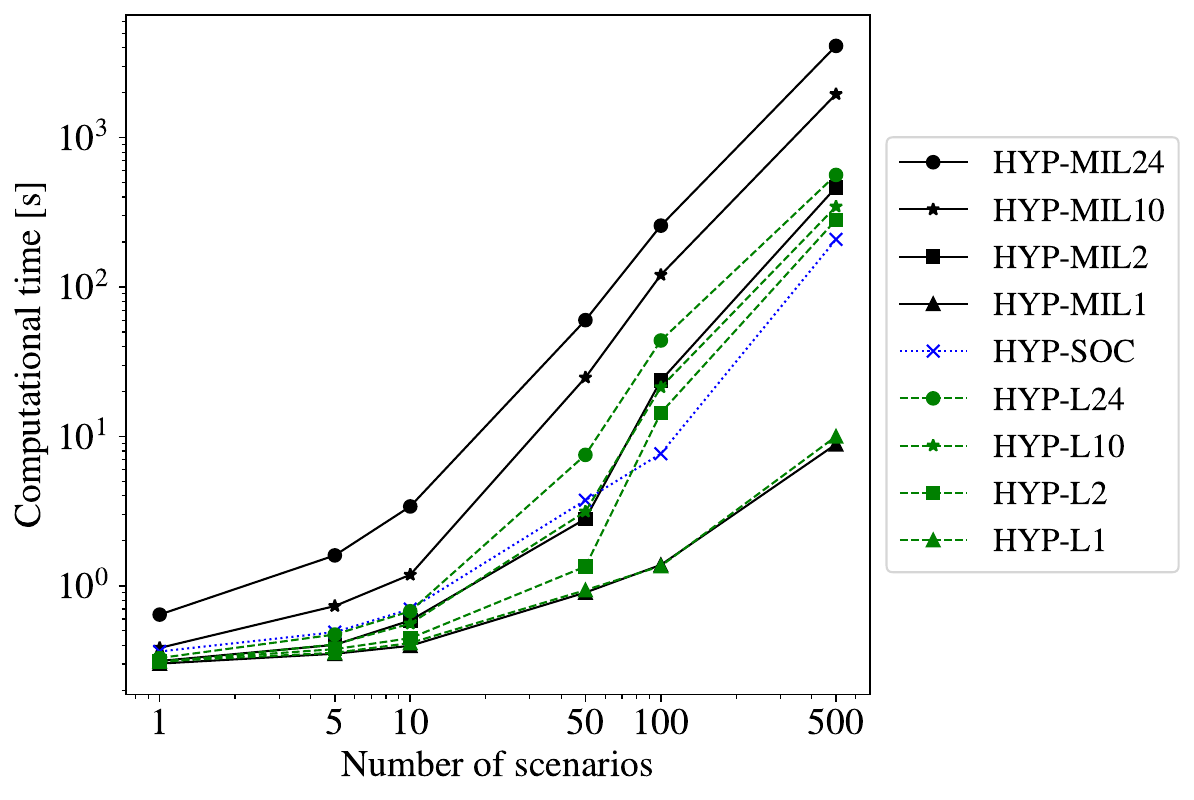}%
\caption{Computational performance of the different models for an increasing number of scenarios. Each data point represents the average computational time obtained from $100$ repetitions, except for $500$ scenarios, where only $10$ repetitions are used due to the increased computational time.} 
\label{fig:comp_time_stochastic}
\end{figure}

Table~\ref{tab:binaries} summarizes the number and type of variables and constraints in each of the hydrogen production curve models.
Figure~\ref{fig:comp_time_stochastic} shows the computational time for the different models in a range of $1$--$500$ scenarios.
\begin{table}[b]
\centering
\caption{Number of variables and constraints for modeling the hydrogen production curve as a function of the number of segments~$\lvert \mathcal{S} \rvert$, time steps $\lvert \mathcal{T} \rvert$, and scenarios $\lvert \Omega \rvert$}
\vspace{0.1cm}
\begin{adjustbox}{width=0.6\textwidth}
\small
\begin{tabular}{lrrrr}
\hline
\multirow{2}{*}{} & \multicolumn{2}{c}{\# Variables} & \multicolumn{2}{c}{\# Constraints} \\ \cline{2-5} 
& Binary        & Continuous       & LP/MILP         & SOC/MISOC        \\ \hline
$\mathrm{HYP}$-$\mathrm{MIL}$ &$\lvert \mathcal{S} \rvert \lvert \mathcal{T} \rvert \lvert \Omega \rvert$               &$\lvert \mathcal{S}\rvert \lvert \mathcal{T} \rvert \lvert \Omega \rvert$                  &$(2 \lvert \mathcal{S}\rvert + 3) \lvert \mathcal{T} \rvert \lvert \Omega \rvert$                 & 0                 \\
$\mathrm{HYP}$-$\mathrm{L}$             &0               &$\lvert \mathcal{T} \rvert \lvert \Omega \rvert$                  &$(\lvert \mathcal{S} \rvert + 3) \lvert \mathcal{T} \rvert \lvert \Omega \rvert $                &0                  \\
$\mathrm{HYP}$-$\mathrm{SOC}$           &0               &$\lvert \mathcal{T} \rvert \lvert \Omega \rvert$                  &$3\lvert \mathcal{T} \rvert \lvert \Omega \rvert$                 &$\lvert \mathcal{T} \rvert \lvert \Omega \rvert$                  \\ \hline
\end{tabular}
\end{adjustbox}
\label{tab:binaries}
\end{table}
All proposed models are solved within $5$ seconds for up to $10$ scenarios. By increasing the number of scenarios, the number of binary variables in the $\mathrm{HYP}$-$\mathrm{MIL}$ model increases comparatively fast, as shown in Table~\ref{tab:binaries}. For $24$ and $10$ segments, the linear relaxation $\mathrm{HYP}$-$\mathrm{L}$ is solved approximately $80$\% faster than the corresponding $\mathrm{HYP}$-$\mathrm{MIL}$, while achieving the same solution. For $100$ scenarios, the $\mathrm{HYP}$-$\mathrm{SOC}$ model is almost three times faster than $\mathrm{HYP}$-$\mathrm{L10}$ and almost two times faster than $\mathrm{HYP}$-$\mathrm{L2}$. Even for $500$ scenarios, the $\mathrm{HYP}$-$\mathrm{SOC}$ model shows a satisfactory computational performance compared to both $\mathrm{HYP}$-$\mathrm{MIL}$ and $\mathrm{HYP}$-$\mathrm{L}$.
This might be attributed to the avoidance of binary variables and the comparably lower number of constraints required to model the hydrogen production curve, even if the constraints are conic instead of linear.
In contrast to model $\mathrm{HYP}$-$\mathrm{SOC}$, $\mathrm{HYP}$-$\mathrm{MISOCP}$, which is presented in Section~OC.1 of the online companion \citep{Online_Companion}, experiences an exponentially increasing computational time that makes it intractable for large applications.

Note that the major part of the computational time is attributed to finding feasible integer solutions for the operational states of the electrolyzer. When fixing the binary  variables to their optimal value, the computational time for model $\mathrm{HYP}$-$\mathrm{SOC}$ including 500 scenarios reduces from 207.9 to 3.7 seconds. This highlights the applicability of the proposed relaxations for real-time rescheduling and the benefit of heuristics finding a set of near-optimal operational states.

For 500 scenarios, the inclusion of linear underestimator~\eqref{eq:DA_linear_underestimator_SOC} in model $\mathrm{HYP}$-$\mathrm{SOC}$ reduces the computational time significantly ($-28\%$), even though it is always inactive at the optimal solution since the relaxation is exact in the considered case study. 
In order to generalize the benefits of adding a linear underestimator to other case studies and model $\mathrm{HYP}$-$\mathrm{L}$, further analyses should be conducted. As this paper focuses on comparing the different modeling approaches in conditions where the relaxations are exact, we leave this aspect for future research.

%% file: 05_Conclusion.tex
\section{Conclusion}
\label{sec:conclusion}
An accurate representation of the nonlinear and nonconvex hydrogen production curve of the electrolyzer, which captures the relationship between power consumption and hydrogen production, is essential for its optimal scheduling. 
The current state-of-the-art modeling approach of the hydrogen production curve is based on piecewise linear approximation. This approach requires carefully selecting the number and the location of linearization points, which impacts the accuracy of dispatch decisions and the computational complexity.
The accuracy increases with the number of linearization segments at the expense of adding one binary variable per segment. For our case study, we found out that using at least ten linearization segments yields sufficiently accurate dispatch decisions. To further highlight computational issues raised by the number of binary variables in the state-of-the-art piecewise linear approximation, we developed a two-stage stochastic program where renewable power production uncertainty is modeled via scenarios. We showed that by growing the number of scenarios and therefore the number of binary variables (indicating active segments for every scenario), the computational time increases significantly, imposing a serious barrier. 

This paper proposes two modeling approaches for the hydrogen production curve based on convex relaxations. The first one, $\mathrm{HYP}$-$\mathrm{L}$, is a linear relaxation of the state-of-the-art piecewise linear approximation $\mathrm{HYP}$-$\mathrm{MIL}$, which does not require binary variables.
Although this leads to a significantly improved computational performance compared to $\mathrm{HYP}$-$\mathrm{MIL}$, it still requires choosing the number and location of linearization points. Additionally, the high number of segments needed to ensure the accuracy of the solution impacts the computational performance of large-scale problems. 

Those barriers are resolved by our second model, $\mathrm{HYP}$-$\mathrm{SOC}$, which is a conic relaxation of a quadratic approximation of the hydrogen production curve. The quadratic approximation can be directly fitted to operational data of the electrolyzer, making it especially suitable in cases with limited information availability on the underlying physics of the hydrogen production curve. Based on a realistic case study, we showed that $\mathrm{HYP}$-$\mathrm{SOC}$ provides a satisfactory trade-off between accuracy of decisions and computational performance for large-scale problems. Between the two proposed relaxations, $\mathrm{HYP}$-$\mathrm{SOC}$ exhibits a slightly better computational performance than $\mathrm{HYP}$-$\mathrm{L}$ for large-scale problems at the expense of moving from a (mixed-integer) linear to a second-order cone problem. We conclude that the linear relaxation provides an appropriate choice if the model type is to be kept linear, otherwise, we suggest using the conic model.

We mathematically proved that the proposed conic relaxation is exact under prevalent operating conditions. Similar proofs can be derived for the linear relaxation. We further presented a heuristic to check the exactness of the conic relaxation a priori based on the wind power availability --- this can be very useful in practice. 
An extreme case was presented to illustrate how a combination of prolonged negative electricity prices, high wind power availability, and a restrictive upper limit for hydrogen production may lead to an inexact solution.

In that case, the feasible region of the problem should be tightened and an exact solution should be restored a posteriori. We were able to show that introducing a linear underestimator to the hydrogen production curve substantially reduces the feasible region. Furthermore, it ensures that an exact solution can always be restored by reducing the power consumption of the electrolyzer. Future research should investigate how the linear underestimator impacts the computational time in cases wherein the relaxation is inexact. As the production of green hydrogen in hours with negative electricity prices constitutes a particularly profitable business case for the electrolyzer, further research should address how an exact yet optimal solution can be restored when other assets or downstream processes are considered.

When disregarding the operational states of the electrolyzer, the two proposed relaxations of the hydrogen production curve are convex. This is a valuable property in cases where global optimality guarantees and meaningful dual variables are desired, e.g., in a market-clearing problem. Future research should explore how the proposed modeling approaches can be integrated into market-oriented applications. Additionally, the conic relaxation and exactness findings may be applicable to other components with similar nonlinear physics, e.g., batteries with a semi-linear linear efficiency curve \citep{Engelhardt2022}.
Finally, future work should focus on accurate modeling of the auxiliary assets and downstream processes, such as the compressor and methanol or ammonia production, which may introduce additional nonconvexities to the optimal scheduling problem of an electrolyzer.

%% file: Appendix_proofs.tex
\section{Mathematical Proofs}
\label{app:proofs}
Let $(\dot{\mathbf{x}}, \dot{\mathbf{y}}, \dot{\mathbf{z}})$ denote a feasible point to problem $\mathrm{HYP}$-$\mathrm{SOC}$. For the following proofs, we assume that the set of binary variables $\mathbf{z}$ associated with the operational states of the electrolyzer is fixed to $\dot{\mathbf{z}}$. In this case, $\mathrm{HYP}$-$\mathrm{SOC}$ reduces to a convex problem. We use $(p^*, h^*, f^*)$ to denote an optimal solution to problem $\mathrm{HYP}$-$\mathrm{SOC}$ for given $\dot{\mathbf{z}}$.

\begin{lemma}\label{Lemma_1}
Suppose the electrolyzer is in standby or off state at time step $t$, such that $\dot{z}_t^{\textrm{sb}} + \dot{z}_t^{\textrm{off}}=1$. Then, the relaxation~\eqref{eq:DA_hydrogen_prod_SOC} is exact.
\end{lemma}

\textit{Proof:} It follows from mutual exclusiveness of the states~\eqref{eq:DA_mutual_states} that $\dot{z}_t^{\textrm{on}} = 0$. Constraints~\eqref{eq:DA_hydrogen_prod_SOC}--\eqref{eq:DA_p_hat_SOC} and \eqref{eq:DA_non_negative_hydrogen} ensure that that the hydrogen production $h_t$ and associated power consumption $\tilde{p}_t$ are equal to zero. Relaxation~\eqref{eq:DA_hydrogen_prod_SOC} is therefore exact when the electrolyzer is in standby or off state.

Without loss of generality, we assume that there exists a time step $\tau \in \mathcal{T}$, for which the electrolyzer is in on-state $\dot{z}_{\tau}^{\textrm{on}}=1$. In this case, it follows from \eqref{eq:DA_power_consumption_tot_SOC} and \eqref{eq:DA_mutual_states} that $p_{\tau} = \tilde{p}_{\tau}$.
In the following, we will prove the exactness of relaxation~\eqref{eq:DA_hydrogen_prod_SOC} under certain assumptions.

\subsection{Proof of Theorem~\ref{Theorem_1}}
\label{sec:proof_theorem_1}
Suppose that $\tau \in \mathcal{H}_{n}$ and $h_{\tau}^* < A {(p_{\tau}^*)}^2 + B p_{\tau}^* + C$, i.e., \eqref{eq:DA_hydrogen_prod_SOC} is inexact. We define
\begin{equation*}
{\Grave{h}_t} =
\begin{cases}
h_{t}^*, & t \neq \tau,\\
{h_{t}^* + \epsilon,} & {t = \tau,}
\end{cases}
\end{equation*}
where $\epsilon$ is a small positive number, such that $\Grave{h}_{\tau} \leq A {(p_{\tau}^*)}^2 + B p_{\tau}^* + C$.
Let the objective function \eqref{eq:DA_objective} be denoted by $U$. As the hydrogen price is positive $\chi > 0$, $U$ is always increasing in $h_t$ for a given $\dot{\mathbf{z}}$. It follows directly that $U(\Grave{h}, f^*) \geq U(h^*, f^*)$. As, by assumption, the total maximum hydrogen production is not constrained by the demand, such that~\eqref{eq:DA_hydrogen_demand} is nonbinding at $h^*$, it is possible to find $\epsilon$ such that ${\Grave{h}}$ is feasible. That contradicts the optimality of $h^*$.
Therefore, relaxation~\eqref{eq:DA_hydrogen_prod_SOC} must be exact at optimum for a given $\dot{\mathbf{z}}$ when the electrolyzer is in on-state.

The exactness of relaxation~\eqref{eq:DA_hydrogen_prod_SOC} for the case when the electrolyzer is in standby or off state follows from Lemma~\eqref{Lemma_1}, which completes the proof.

\subsection{Proof of Theorem~\ref{Theorem_2}}
\label{sec:proof_theorem_2}
It follows directly from the power balance of the hybrid power plant~\eqref{eq:DA_balance}, that $f^* = W - p^*$.
Suppose that in sub-period $n$, constraints~\eqref{eq:DA_hydrogen_prod_SOC} are binding at optimum for all time steps, such that $h_t^* = A (p_t^*)^2 + B p_t^* + C, \forall t \in \mathcal{H}_n$, and that $\tau \in \mathcal{H}_n$.
We define
\begin{equation*}
(\Acute{p}_t, \Acute{f}_t) =
\begin{cases}
(p_{t}^*, f_{t}^*) , & t \neq \tau,\\
{(p_{t}^* + \epsilon, f_{t}^* - \epsilon ),} & {t = \tau,}
\end{cases}
\end{equation*}
where $\epsilon$ is a small positive number. Suppose $(\Acute{p}, h^*, \Acute{f})$ denotes a feasible point to $\mathrm{HYP}$-$\mathrm{SOC}$, such that $\Acute{f} = W - \Acute{p}$ and $\Acute{f}_{\tau} < f_{\tau}^*$. Then $h_{\tau}^* < A (\Acute{p}_{\tau})^2 + B \Acute{p}_{\tau} + C$, i.e., relaxation~\eqref{eq:DA_hydrogen_prod_SOC} is inexact. 
As the objective function $U$ is increasing in $f$ when $\lambda > 0$, it follows that $U(h^*,\Acute{f}) < U(h^*,f^*)$, which contradicts optimality of $\Acute{f}$. This proves the optimality of $(p^*, h^*, f^*)$ and the exactness of the relaxation~\eqref{eq:DA_hydrogen_prod_SOC} for a given $\dot{\mathbf{z}}$.

The exactness of relaxation~\eqref{eq:DA_hydrogen_prod_SOC} for the case when the electrolyzer is in standby or off state follows from Lemma~\eqref{Lemma_1}, which completes the proof.

\subsection{Proof of Theorem~\ref{Theorem_3}}
\label{sec:proof_theorem_3}
Suppose in sub-period $n$ the set of nonpositive power prices is nonempty $\mathcal{T}_{n}^- \neq \emptyset$, and 
the relaxation~\eqref{eq:DA_hydrogen_prod_SOC} is binding at optimum for all time steps with nonpositive prices, such that $h_t^* = A (p_t^*)^2 + B p_t^* + C, \forall t \in \mathcal{T}_{n}^-$. We define $(\Acute{p}_t, \Acute{f}_t)$ as in the proof of Theorem~\ref{Theorem_2} for $\tau \in \mathcal{T}_{n}^-$, such that $h_{\tau}^* < A (\Acute{p}_{\tau})^2 + B \Acute{p}_{\tau} + C$, i.e., constraint~\eqref{eq:DA_hydrogen_prod_SOC} is nonbinding. As the objective function $U$ is decreasing in $f$ when $\lambda < 0$, it follows that $U(\Grave{f}) \geq U(f^*)$, which contradicts the optimality of $f^*$. We conclude that in sub-period $n$, for a given $\dot{\mathbf{z}}$, there exists at least one time step $\tau \in \mathcal{T}_{n}^-$, for which the optimal solution $(p^*, h^*, f^*)$ has a nonzero feasibility gap, i.e., relaxation $\mathrm{HYP}$-$\mathrm{SOC}$ is inexact.

%% file: Acknowledgment.tex
\section*{Acknowledgement}
This work is supported in part by the Energy Cluster Denmark through the ``Sustainable P2X Business Model" project.
We are grateful to Alexandra Lüth (Copenhagen Business School) for initial discussions on the conic relaxation of the hydrogen production curve and to Josh Taylor (New Jersey Institute of Technology) for suggesting the use of the linear underestimator.
We would like to thank Jens Jakob Sørensen, Alexander Holm Kiilerich, Mathias Stolpe (all Ørsted), and Roar Hestbek Nicolaisen (Hybrid Greentech) for collaborations, thoughtful discussions, and constructive feedback. We thank Anubhav Ratha (Vestas) for discussions on the exactness proofs and providing feedback on the manuscript. Our final gratitude goes to two anonymous reviewers for providing critical feedback that helped us further improve our work.